\documentclass[12pt,leqno,bezier]{amsart}
\usepackage{latexsym,amssymb,amscd,amsthm,amstext,amsfonts}
\unitlength=1mm\linethickness{0.8pt}
\def\emline#1#2#3#4#5#6{\put(#1,#2){\special{em:moveto}}%
       \put(#4,#5){\special{em:lineto}}}\def\newpic#1{}

\def\guy{\special{em:linewidth 0.5pt}\linethickness{0.5pt}
\begin{picture}(3,3)\emline{1}{0}{1}{1}{2}{2}
\emline{1}{2}{3}{0}{3}{4}\emline{1}{2}{5}{2}{3}{6}\end{picture}}
\def\guc{
\begin{picture}(3,3)\emline{1}{3}{1}{1}{1}{2}
\emline{1}{1}{3}{0}{0}{4}\emline{1}{1}{5}{2}{0}{6}\end{picture}}
\def\guo{\special{em:linewidth 0.4pt}\linethickness{0.4pt}
\begin{picture}(2,4)\put(1,2){\circle{2}}
\emline{1}{1}{1}{1}{0}{2}\emline{1}{3}{3}{1}{4}{4}\end{picture}}
\def\guv{\special{em:linewidth 0.4pt}\linethickness{0.4pt}
\begin{picture}(7,2)\put(1,0.5){\vector(-1,0){0.1}}
\emline{5}{0.5}{1}{1}{0.5}{2}\put(1,2.5){\vector(-1,0){0.1}}
\emline{5}{2.5}{3}{1}{2.5}{4}\end{picture}}
\def\gua{\special{em:linewidth 0.4pt}\linethickness{0.4pt}
\begin{picture}(5,6)\put(4,4){\oval(2,2)[t]}\put(2,4){\oval(2,2)[b]}
\put(3.50,3){\oval(3,4)[b]}\emline{3}{0}{1}{3}{1}{2}
\emline{1}{6}{3}{1}{4}{4}\emline{4}{6}{5}{4}{5}{6}
\emline{5}{4}{7}{5}{3}{8}\end{picture}}
\def\gur{\special{em:linewidth 0.4pt}\linethickness{0.4pt}
\begin{picture}(3,6)\emline{3}{6}{1}{3}{0}{2}
\emline{1}{6}{3}{1}{3}{4}\put(1,3){\circle*{1.33}}\end{picture}}
\def\guj{\begin{picture}(4,4)\emline{2}{0}{1}{2}{4}{2}
\put(2,2){\circle*{1}}\end{picture}}
\def\gub{\special{em:linewidth 0.4pt}\unitlength 1mm\linethickness{0.4pt}
\begin{picture}(5,3)\emline{1}{3}{1}{4}{0}{2}
\emline{4}{3}{3}{3}{2}{4}\emline{2}{1}{5}{1}{0}{6}\end{picture}\,}
\def\gcr{\special{em:linewidth 0.4pt}\unitlength 1mm\linethickness{0.4pt}
\begin{picture}(2,3)\emline{1}{0}{1}{1}{2}{2}
\put(1,2){\circle*{1}}\end{picture}}

\def\guA{\special{em:linewidth 0.4pt}\unitlength 1mm\linethickness{0.4pt}
\begin{center}\begin{picture}(66,20)
\put(51,12){\circle{6}}
\emline{51}{14.67}{1}{48}{11.67}{2}
\emline{53}{14}{3}{48.67}{9.67}{4}
\emline{53.67}{12}{5}{50.67}{9}{6}
\emline{46}{19}{7}{46}{17}{8}
\emline{46}{17}{9}{48.33}{14.67}{10}
\emline{48}{19}{11}{50.33}{14.67}{12}
\emline{54}{19}{13}{52}{15}{14}
\emline{56}{19}{15}{56}{17}{16}
\emline{56}{17}{17}{53}{14}{18}
\emline{46}{5}{19}{46}{7}{20}
\emline{46}{7}{21}{49}{10}{22}
\emline{48}{5}{23}{50}{9.33}{24}
\emline{56}{5}{25}{56}{7}{26}
\emline{56}{7}{27}{53}{10}{28}
\emline{54}{5}{29}{52.33}{9}{30}
\put(51,19){\makebox(0,0)[cc]{$\cdots$}}
\put(51,5){\makebox(0,0)[cc]{$\cdots$}}
\put(61,12){\makebox(0,0)[cc]{$\in$}}
\put(66,12){\makebox(0,0)[lc]{$m\mapsto n.$}}
\put(51,20){\makebox(0,0)[cb]{$\stackrel{m}{\overbrace{\qquad\quad\;}}$}}
\put(51,4){\makebox(0,0)[ct]{$\stackrel{\underbrace{\qquad\quad\;}}{n}$}}
\put(15,12){\makebox(0,0)[rc]{\begin{tabular}{c}\textsc{Figure} 1.
\end{tabular}}}\end{picture}\end{center}}

\def\guB{$$\begin{picture}(120,22) 
\emline{5}{20}{1}{5}{15}{2}
\emline{19}{10}{3}{19}{5}{4}
\emline{33}{20}{5}{33}{15}{6}
\emline{33}{10}{7}{33}{5}{8}
\emline{69}{20}{9}{69}{15}{10}
\emline{69}{10}{11}{69}{5}{12}
\emline{46}{20}{13}{46}{5}{14}
\emline{59}{5}{15}{59}{20}{16}
\emline{84}{20}{17}{84}{5}{18}
\emline{94}{5}{19}{94}{20}{20}
\emline{109}{20}{21}{109}{5}{22}
\emline{119}{5}{23}{119}{20}{24}
\put(119,13){\circle*{2}}
\put(109,15){\circle*{2}}
\put(109,10){\circle*{2}}
\put(84,10){\circle*{2}}
\put(84,15){\circle*{2}}
\put(69,15){\circle*{2}}
\put(69,10){\circle*{2}}
\put(59,10){\circle*{2}}
\put(59,15){\circle*{2}}
\put(33,15){\circle*{2}}
\put(33,10){\circle*{2}}
\put(19,10){\circle*{2}}
\put(5,15){\circle*{2}}
\put(64,13){\makebox(0,0)[cc]{\large$\approx$}}
\put(89,13){\makebox(0,0)[cc]{\large$\approx$}}
\put(114,13){\makebox(0,0)[cc]{\large$\approx$}}
\put(5,22){\makebox(0,0)[cb]{killer}}
\put(19,22){\makebox(0,0)[cb]{creator}}
\put(33,22){\makebox(0,0)[cb]{const}}
\put(46,22){\makebox(0,0)[cb]{id}}
\put(64,22){\makebox(0,0)[cb]{\textit{nil}potent}}
\put(89,22){\makebox(0,0)[cb]{\textit{uni}potent}}
\put(114,22){\makebox(0,0)[cb]{\textit{idem}potent}}
\put(5,4){\makebox(0,0)[ct]{$1\mapsto 0$}}
\put(19,4){\makebox(0,0)[ct]{$0\mapsto 1$}}
\put(33,4){\makebox(0,0)[ct]{$1\mapsto 1$}}
\put(113.67,4){\makebox(0,0)[ct]{$1\mapsto 1$}}
\end{picture}$$$$\begin{picture}(46,10) 
\put(7,15){\oval(4,10)[b]}
\put(42,5){\oval(4,10)[t]}
\emline{42}{10}{1}{42}{15}{2}
\emline{7}{5}{3}{7}{10}{4}
\put(7,4){\makebox(0,0)[ct]{$2\mapsto 1$}}
\put(42,4){\makebox(0,0)[ct]{$1\mapsto 2$}}
\put(23,0){\makebox(0,0)[ct]{\begin{tabular}{c}\textsc{Figure} 2.
Some plants and some axioms.\end{tabular}}}\end{picture}$$\bigskip}

\def\guC{$$\begin{picture}(114,22)   
\put(7,18){\oval(4,8)[b]}
\put(13,10){\oval(4,8)[t]}
\put(30,18){\oval(4,6)[b]}
\put(30,10){\oval(4,6)[t]}
\emline{30}{13}{1}{30}{15}{2}
\emline{13}{18}{3}{13}{14}{4}
\emline{7}{14}{5}{7}{10}{6}
\put(47,10){\oval(4,8)[t]}
\emline{47}{14}{7}{47}{18}{8}
\put(53,18){\oval(4,8)[b]}
\emline{53}{14}{9}{53}{10}{10}
\put(70,14){\oval(4,4)[t]}
\put(74,14){\oval(4,4)[b]}
\emline{76}{14}{11}{76}{18}{12}
\emline{70}{18}{13}{70}{16}{14}
\emline{68}{14}{15}{68}{10}{16}
\emline{74}{10}{17}{74}{12}{18}
\put(91,14){\circle{4}}
\emline{91}{18}{19}{91}{16}{20}
\emline{91}{12}{21}{91}{10}{22}
\put(108,14){\oval(4,4)[b]}
\put(112,14){\oval(4,4)[t]}
\emline{112}{16}{23}{112}{18}{24}
\emline{106}{18}{25}{106}{14}{26}
\emline{108}{12}{27}{108}{10}{28}
\emline{114}{10}{29}{114}{14}{30}
\put(10,9){\makebox(0,0)[ct]{$3\mapsto 3$}}
\put(30,9){\makebox(0,0)[ct]{$2\rightarrow 2$}}
\put(50,9){\makebox(0,0)[ct]{$3\mapsto 3$}}
\put(72,9){\makebox(0,0)[ct]{$2\mapsto 2$}}
\put(91,9){\makebox(0,0)[ct]{$1\mapsto 1$}}
\put(110,9){\makebox(0,0)[ct]{$2\mapsto 2$}}
\put(57,3){\makebox(0,0)[ct]{\begin{tabular}{c}\textsc{Figure} 3.
Two-letters derived plants in a garden.\end{tabular}}}\end{picture}$$
\smallskip}

\def\guD{$$\begin{picture}(84,22)        
\emline{1}{10}{1}{1}{18}{2}
\emline{26}{18}{3}{26}{14}{4}
\put(26,14){\circle*{2}}
\put(53,18){\oval(4,8)[b]}
\emline{53}{14}{5}{53}{10}{6}
\put(82,10){\oval(4,8)[t]}
\emline{82}{14}{7}{82}{18}{8}
\put(1,9){\makebox(0,0)[ct]{id}}
\put(26,9){\makebox(0,0)[ct]{$1\mapsto 0$}}
\put(53,9){\makebox(0,0)[ct]{$2\mapsto 1$}}
\put(82,9){\makebox(0,0)[ct]{$1\mapsto 2$}}
\put(42,3){\makebox(0,0)[ct]{\begin{tabular}{c}\textsc{Figure} 4.
Alphabet of plants.\end{tabular}}}\end{picture}$$}

\def\guE{\special{em:linewidth 0.4pt}\unitlength 1mm\linethickness{0.4pt}
$$\begin{picture}(46,44)   
\put(14,42){\circle{2}}
\put(12.50,41){\oval(3,4)[b]}
\emline{14}{43}{1}{14}{44}{2}
\emline{11}{44}{3}{11}{41}{4}
\emline{13}{38}{5}{13}{39}{6}
\put(34,42){\oval(2,2)[t]}
\put(32,42){\oval(2,2)[b]}
\put(33.50,41){\oval(3,4)[b]}
\emline{35}{41}{7}{35}{42}{8}
\emline{31}{44}{9}{31}{42}{10}
\emline{34}{44}{11}{34}{43}{12}
\emline{33.33}{38}{13}{33.33}{39}{14}
\put(24,27){\circle{2}}
\emline{24}{30}{15}{24}{28}{16}
\emline{24}{26}{17}{24}{24}{18}
\emline{21}{30}{19}{21}{27}{20}
\put(21,27){\circle*{1}}
\put(3,30){\oval(4,6)[b]}
\emline{3}{24}{21}{3}{27}{22}
\put(44,27){\oval(2,4)[t]}
\put(43,27){\oval(4,4)[b]}
\put(43,27){\circle*{1}}
\emline{41}{30}{23}{41}{27}{24}
\emline{44}{30}{25}{44}{29}{26}
\emline{43}{24}{27}{43}{25}{28}
\put(34,13){\oval(2,4)[t]}
\emline{34}{16}{29}{34}{15}{30}
\emline{31}{16}{31}{31}{13}{32}
\put(31,13){\circle*{1}}
\put(33,13){\circle*{1}}
\emline{35}{10}{33}{35}{13}{34}
\emline{14}{10}{35}{14}{16}{36}
\emline{12}{16}{37}{12}{13}{38}
\put(12,13){\circle*{1}}
\put(23,3){\makebox(0,0)[ct]{\begin{tabular}{c}\textsc{Figure} 5.
Small fragment of a garden of type $2\mapsto 1.$
\end{tabular}}}\end{picture}$$}

\def\guF{$$\begin{picture}(116,13)    
\emline{111}{1}{1}{111}{9}{2}
\put(26,5){\circle{4}}
\emline{26}{7}{3}{26}{9}{4}
\emline{26}{1}{5}{26}{3}{6}
\emline{71}{1}{7}{71}{9}{8}
\emline{46}{1}{9}{46}{9}{10}
\put(91,9){\oval(4,8)[b]}
\emline{91}{1}{11}{91}{5}{12}
\emline{68}{9}{13}{68}{5}{14}
\put(68,5){\circle*{1.33}}
\emline{114}{9}{15}{114}{5}{16}
\put(114,5){\circle*{1.33}}
\put(41,6){\vector(4,-3){0.2}}
\emline{33}{6}{17}{34.82}{7.40}{18}
\emline{34.82}{7.40}{19}{36.64}{7.98}{20}
\emline{36.64}{7.98}{21}{38.45}{7.74}{22}
\emline{38.45}{7.74}{23}{41}{6}{24}
\put(33,4){\vector(-4,3){0.2}}
\emline{41}{4}{25}{39.18}{2.60}{26}
\emline{39.18}{2.60}{27}{37.36}{2.02}{28}
\emline{37.36}{2.02}{29}{35.55}{2.26}{30}
\emline{35.55}{2.26}{31}{33}{4}{32}
\put(76,6){\vector(-4,-3){0.2}}
\emline{84}{6}{33}{82.18}{7.40}{34}
\emline{82.18}{7.40}{35}{80.36}{7.98}{36}
\emline{80.36}{7.98}{37}{78.55}{7.74}{38}
\emline{78.55}{7.74}{39}{76}{6}{40}
\put(84,4){\vector(4,3){0.2}}
\emline{76}{4}{41}{77.82}{2.60}{42}
\emline{77.82}{2.60}{43}{79.64}{2.02}{44}
\emline{79.64}{2.02}{45}{81.45}{2.26}{46}
\emline{81.45}{2.26}{47}{84}{4}{48}
\put(106,6){\vector(4,-3){0.2}}
\emline{98}{6}{49}{99.82}{7.40}{50}
\emline{99.82}{7.40}{51}{101.64}{7.98}{52}
\emline{101.64}{7.98}{53}{103.45}{7.74}{54}
\emline{103.45}{7.74}{55}{106}{6}{56}
\put(98,4){\vector(-4,3){0.2}}
\emline{106}{4}{57}{104.18}{2.60}{58}
\emline{104.18}{2.60}{59}{102.36}{2.02}{60}
\emline{102.36}{2.02}{61}{100.55}{2.26}{62}
\emline{100.55}{2.26}{63}{98}{4}{64}
\put(15,9){\makebox(0,0)[cb]{$u\times\id$}}
\put(15,1){\makebox(0,0)[ct]{$\varepsilon\times\id$}}
\put(37,9){\makebox(0,0)[cb]{$m$}}
\put(37,1){\makebox(0,0)[ct]{$\triangle$}}
\put(80,1){\makebox(0,0)[ct]{$L$}}
\put(80,9){\makebox(0,0)[cb]{$l$}}
\put(102,9){\makebox(0,0)[cb]{$r$}}
\put(102,1){\makebox(0,0)[ct]{$R$}}
\put(4,4){\oval(4,4)[b]}
\put(4,6){\oval(4,4)[t]}
\emline{6}{6}{65}{6}{4}{66}
\emline{4}{9}{67}{4}{8}{68}
\emline{4}{1}{69}{4}{2}{70}
\put(2,4){\circle*{1}}
\put(2,6){\circle*{1}}
\put(19,6){\vector(4,-3){0.2}}
\emline{11}{6}{71}{12.82}{7.40}{72}
\emline{12.82}{7.40}{73}{14.64}{7.98}{74}
\emline{14.64}{7.98}{75}{16.45}{7.74}{76}
\emline{16.45}{7.74}{77}{19}{6}{78}
\put(11,4){\vector(-4,3){0.2}}
\emline{19}{4}{79}{17.18}{2.60}{80}
\emline{17.18}{2.60}{81}{15.36}{2.02}{82}
\emline{15.36}{2.02}{83}{13.55}{2.26}{84}
\emline{13.55}{2.26}{85}{11}{4}{86}
\end{picture}$$}

\def\guG{\special{em:linewidth 0.4pt}\unitlength 1mm\linethickness{0.4pt}
$$\begin{picture}(110,18)     
\put(40,11){\oval(4,4)[t]}
\put(38,11){\oval(8,8)[b]}
\emline{34}{15}{1}{34}{11}{2}
\emline{40}{15}{3}{40}{13}{4}
\emline{38}{5}{5}{38}{7}{6}
\put(38,11){\circle*{1.33}}
\put(56,15){\oval(4,10)[b]}
\emline{56}{5}{7}{56}{10}{8}
\put(88,11){\circle{4}}
\put(86,9){\oval(4,4)[b]}
\emline{88}{13}{9}{88}{15}{10}
\emline{84}{15}{11}{84}{9}{12}
\emline{86}{5}{13}{86}{7}{14}
\put(104,11){\oval(4,4)[b]}
\put(108,11){\oval(4,4)[t]}
\put(107,9){\oval(6,6)[b]}
\emline{107}{4}{15}{107}{6}{16}
\emline{110}{9}{17}{110}{11}{18}
\emline{108}{15}{19}{108}{13}{20}
\emline{102}{15}{21}{102}{11}{22}
\emline{18}{15}{23}{18}{4}{24}
\emline{4}{15}{25}{4}{12}{26}
\put(4,7){\oval(4,10)[t]}
\emline{6}{4}{27}{6}{7}{28}
\put(2,7){\circle*{1.33}}
\put(26,10){\makebox(0,0)[cc]{$\Longrightarrow$}}
\put(46,16){\makebox(0,0)[cb]{$\stackrel{\hbox{expanded c}}{\overbrace
{\qquad\qquad\qquad\qquad}}$}}
\put(16,10){\vector(1,0){0.2}}
\emline{8}{10}{29}{16}{10}{30}
\put(52,10){\vector(1,0){0.2}}
\emline{44}{10}{31}{52}{10}{32}
\put(92,10){\vector(-1,0){0.2}}
\emline{100}{10}{33}{92}{10}{34}
\put(12,11){\makebox(0,0)[cb]{c}}
\put(48,11){\makebox(0,0)[cb]{id $\times$ c}}
\put(96,11){\makebox(0,0)[cb]{a}}
\put(59,10){\makebox(0,0)[lc]{;}}
\put(55,3){\makebox(0,0)[ct]{\begin{tabular}{c}\textsc{Figure} 6. Two
more footpaths needed.\end{tabular}}}\end{picture}$$}

\def\guH{\special{em:linewidth 0.4pt}\unitlength 1mm\linethickness{0.4pt}
$$\begin{picture}(40,22)         
\put(20,20){\circle{2}}
\emline{20}{21}{1}{20}{22}{2}
\emline{20}{19}{3}{20}{18}{4}
\emline{20}{10}{5}{20}{6}{6}
\put(3,7){\oval(2,4)[t]}
\put(37,7){\oval(2,4)[t]}
\emline{37}{9}{7}{37}{10}{8}
\emline{36}{7}{9}{36}{6}{10}
\put(38,7){\circle*{1.33}}
\emline{3}{10}{11}{3}{9}{12}
\emline{4}{6}{13}{4}{7}{14}
\put(2,7){\circle*{1.33}}
\put(12,8){\makebox(0,0)[cc]{\large$\stackrel{c}{\rightarrow}$}}
\put(3,13){\vector(-1,-4){0.2}}
\emline{15}{21}{15}{10.50}{21}{16}
\emline{10.50}{21}{17}{8.63}{20.50}{18}
\emline{8.63}{20.50}{19}{7}{19.67}{20}
\emline{7}{19.67}{21}{5.63}{18.50}{22}
\emline{5.63}{18.50}{23}{4.50}{17}{24}
\emline{4.50}{17}{25}{3.63}{15.17}{26}
\emline{3.63}{15.17}{27}{3}{13}{28}
\put(37,13){\vector(1,-4){0.2}}
\emline{25}{21}{29}{27.55}{21.17}{30}
\emline{27.55}{21.17}{31}{29.79}{21.03}{32}
\emline{29.79}{21.03}{33}{31.73}{20.59}{34}
\emline{31.73}{20.59}{35}{33.37}{19.85}{36}
\emline{33.37}{19.85}{37}{34.71}{18.81}{38}
\emline{34.71}{18.81}{39}{35.75}{17.46}{40}
\emline{35.75}{17.46}{41}{36.48}{15.81}{42}
\emline{36.48}{15.81}{43}{37}{13}{44}
\put(19,12){\vector(1,-2){0.2}}
\emline{19}{16}{45}{18.03}{14.33}{46}
\emline{18.03}{14.33}{47}{19}{12}{48}
\put(21,16){\vector(-1,2){0.2}}
\emline{21}{12}{49}{21.97}{13.67}{50}
\emline{21.97}{13.67}{51}{21}{16}{52}
\put(15,19){\vector(1,0){0.2}}
\emline{6}{13}{53}{6.58}{15.05}{54}
\emline{6.58}{15.05}{55}{7.61}{16.67}{56}
\emline{7.61}{16.67}{57}{9.08}{17.87}{58}
\emline{9.08}{17.87}{59}{11}{18.65}{60}
\emline{11}{18.65}{61}{15}{19}{62}
\put(25,19){\vector(-1,0){0.2}}
\emline{34}{13}{63}{33.32}{15.05}{64}
\emline{33.32}{15.05}{65}{32.24}{16.67}{66}
\emline{32.24}{16.67}{67}{30.75}{17.87}{68}
\emline{30.75}{17.87}{69}{28.86}{18.65}{70}
\emline{28.86}{18.65}{71}{25}{19}{72}
\put(17,14){\makebox(0,0)[rc]{$m$}}
\put(23,14){\makebox(0,0)[lc]{$\triangle$}}
\put(4,17){\makebox(0,0)[rb]{$l$}}
\put(37,17){\makebox(0,0)[lb]{$r$}}
\put(8,16){\makebox(0,0)[lt]{$L$}}
\put(32,16){\makebox(0,0)[rt]{$R$}}
\put(20,4){\makebox(0,0)[ct]{\begin{tabular}{c}\textsc{Figure} 7.
Footpaths `$m$' and `$c\circ l$' almost coincide on a plant $\guo$ ,\\
$c\circ l\approx m$\quad\&\quad $L\approx\triangle\circ c.$
\end{tabular}}}\end{picture}$$}

\def\guI{$$\begin{picture}(54,17)       
\put(5,8.17){\oval(6,12.04)[t]}
\put(51,8.17){\oval(6,12.04)[t]}
\emline{5}{17}{1}{5}{14.19}{2}
\emline{28}{17}{3}{28}{5.16}{4}
\emline{51}{17}{5}{51}{14.19}{6}
\emline{2}{8.17}{7}{2}{5.16}{8}
\emline{54}{5.16}{9}{54}{8.17}{10}
\put(48,8.17){\circle*{1.33}}
\put(8,8.17){\circle*{1.33}}
\put(18,11.18){\makebox(0,0)[cc]{\large$\approx$}}
\put(38,11.18){\makebox(0,0)[cc]{\large$\approx$}}
\put(27,3){\makebox(0,0)[ct]{\begin{tabular}{c}\textsc{Figure} 8.
Strict mitosis.\end{tabular}}}\end{picture}$$}

\def\guJ{
$$\begin{picture}(64,26)           
\emline{5}{18}{1}{8}{15}{2}
\emline{8}{15}{3}{8}{14}{4}
\emline{11}{18}{5}{8}{15}{6}
\emline{7}{18}{7}{6}{17}{8}
\emline{9}{18}{9}{7}{16}{10}
\emline{30}{9}{11}{33}{6}{12}
\emline{33}{6}{13}{33}{5}{14}
\emline{36}{9}{15}{33}{6}{16}
\emline{32}{9}{17}{31}{8}{18}
\emline{34}{9}{19}{35}{8}{20}
\emline{55}{18}{21}{58}{15}{22}
\emline{58}{15}{23}{58}{14}{24}
\emline{61}{18}{25}{58}{15}{26}
\emline{59}{18}{27}{60}{17}{28}
\emline{57}{18}{29}{59}{16}{30}
\emline{22}{29}{31}{25}{26}{32}
\emline{25}{26}{33}{25}{25}{34}
\emline{26}{29}{35}{24}{27}{36}
\emline{28}{29}{37}{25}{26}{38}
\emline{24}{29}{39}{25}{28}{40}
\emline{38}{29}{41}{41}{26}{42}
\emline{41}{26}{43}{41}{25}{44}
\emline{44}{29}{45}{41}{26}{46}
\emline{40}{29}{47}{42}{27}{48}
\emline{42}{29}{49}{41}{28}{50}
\put(36,27){\vector(1,0){0.2}}
\emline{30}{27}{51}{36}{27}{52}
\put(20,25){\vector(1,1){0.2}}
\emline{13}{18}{53}{20}{25}{54}
\put(53,18){\vector(1,-1){0.2}}
\emline{46}{25}{55}{53}{18}{56}
\put(27,7){\vector(2,-1){0.2}}
\emline{13}{14}{57}{27}{7}{58}
\put(53,14){\vector(2,1){0.2}}
\emline{39}{7}{59}{53}{14}{60}
\put(20,11){\makebox(0,0)[lb]{$a$}}
\put(46,11){\makebox(0,0)[rb]{$a$}}
\put(17,21){\makebox(0,0)[lt]{$a\times\id$}}
\put(50,22){\makebox(0,0)[lb]{$\id\times a$}}
\put(33,28){\makebox(0,0)[cb]{$a$}}
\put(32,2){\makebox(0,0)[ct]{\textsc{Figure} 9. The MacLane \& Stasheff
pentagon.}}\end{picture}$$}

\def\guK{\special{em:linewidth 0.4pt}\unitlength 1mm\linethickness{0.4pt}
\begin{center}\begin{picture}(100,35)         
\put(1,5){\framebox(4,20)[cc]{}}
\put(24,5){\framebox(8,20)[cc]{}}
\put(51,5){\framebox(14,20)[cc]{}}
\put(88,5){\framebox(12,20)[cc]{}}
\put(94,15){\circle{4}}
\put(61,15){\circle{4}}
\put(55,15){\circle{4}}
\put(28,15){\circle{4}}
\put(58,13){\oval(6,8)[b]}
\put(58,17){\oval(6,8)[t]}
\put(95,11){\oval(6,4)[b]}
\put(92,13){\oval(4,4)[b]}
\put(96,17){\oval(4,4)[t]}
\put(93,19){\oval(6,4)[t]}
\emline{93}{23}{1}{93}{21}{2}
\emline{95}{7}{3}{95}{9}{4}
\emline{98}{11}{5}{98}{17}{6}
\emline{90}{13}{7}{90}{19}{8}
\emline{58}{23}{9}{58}{21}{10}
\emline{58}{7}{11}{58}{9}{12}
\emline{28}{7}{13}{28}{13}{14}
\emline{28}{17}{15}{28}{23}{16}
\emline{3}{23}{17}{3}{7}{18}
\put(30,27){\vector(1,-4){0.2}}
\emline{26}{27}{19}{26.63}{29.11}{20}
\emline{26.63}{29.11}{21}{27.25}{30.44}{22}
\emline{27.25}{30.44}{23}{27.88}{30.98}{24}
\emline{27.88}{30.98}{25}{28.50}{30.75}{26}
\emline{28.50}{30.75}{27}{29.13}{29.73}{28}
\emline{29.13}{29.73}{29}{30}{27}{30}
\put(97,27){\vector(1,-2){0.2}}
\emline{90}{27}{31}{90.99}{29.01}{32}
\emline{90.99}{29.01}{33}{91.99}{30.32}{34}
\emline{91.99}{30.32}{35}{93}{30.94}{36}
\emline{93}{30.94}{37}{94.04}{30.88}{38}
\emline{94.04}{30.88}{39}{95.08}{30.11}{40}
\emline{95.08}{30.11}{41}{97}{27}{42}
\put(19,16){\vector(2,-1){0.2}}
\emline{10}{16}{43}{11.82}{17.32}{44}
\emline{11.82}{17.32}{45}{13.67}{17.94}{46}
\emline{13.67}{17.94}{47}{15.55}{17.88}{48}
\emline{15.55}{17.88}{49}{19}{16}{50}
\put(10,14){\vector(-2,1){0.2}}
\emline{19}{14}{51}{17.07}{12.68}{52}
\emline{17.07}{12.68}{53}{15.17}{12.06}{54}
\emline{15.17}{12.06}{55}{13.30}{12.13}{56}
\emline{13.30}{12.13}{57}{10}{14}{58}
\put(46,16){\vector(2,-1){0.2}}
\emline{37}{16}{59}{38.82}{17.32}{60}
\emline{38.82}{17.32}{61}{40.67}{17.94}{62}
\emline{40.67}{17.94}{63}{42.55}{17.88}{64}
\emline{42.55}{17.88}{65}{46}{16}{66}
\put(37,14){\vector(-2,1){0.2}}
\emline{46}{14}{67}{44.07}{12.68}{68}
\emline{44.07}{12.68}{69}{42.17}{12.06}{70}
\emline{42.17}{12.06}{71}{40.30}{12.13}{72}
\emline{40.30}{12.13}{73}{37}{14}{74}
\put(81,16){\vector(2,-1){0.2}}
\emline{72}{16}{75}{73.82}{17.32}{76}
\emline{73.82}{17.32}{77}{75.67}{17.94}{78}
\emline{75.67}{17.94}{79}{77.55}{17.88}{80}
\emline{77.55}{17.88}{81}{81}{16}{82}
\put(72,14){\vector(-2,1){0.2}}
\emline{81}{14}{83}{79.07}{12.68}{84}
\emline{79.07}{12.68}{85}{77.17}{12.06}{86}
\emline{77.17}{12.06}{87}{75.30}{12.13}{88}
\emline{75.30}{12.13}{89}{72}{14}{90}
\put(30,29){\makebox(0,0)[lc]{$\sigma$}}
\put(97,29){\makebox(0,0)[lc]{$\sigma$}}
\put(15,19){\makebox(0,0)[cb]{$\triangle$}}
\put(14,11){\makebox(0,0)[ct]{$m$}}
\put(41,11){\makebox(0,0)[ct]{$m\times m$}}
\put(76,11){\makebox(0,0)[ct]{$a(a^{-1}\times\id)$}}
\put(77,19){\makebox(0,0)[cb]{$(a\times\id)a^{-1}$}}
\put(42,19){\makebox(0,0)[cb]{$\triangle\times\triangle$}}
\put(50,3){\makebox(0,0)[ct]{\begin{tabular}{c}\textsc{Figure} 10.
A $\sigma$-braided bigebra.\end{tabular}}}\end{picture}\end{center}}

\def\guL{\special{em:linewidth 0.4pt}\unitlength 2mm\linethickness{0.4pt}
$$\begin{picture}(46,40)         
\put(14,38.65){\circle{2}}
\put(12.50,37.60){\oval(3,4.20)[b]}
\emline{14}{39.70}{1}{14}{40.75}{2}
\emline{11}{40.75}{3}{11}{37.60}{4}
\emline{13}{34.45}{5}{13}{35.50}{6}
\put(34,38.65){\oval(2,2.10)[t]}
\put(32,38.65){\oval(2,2.10)[b]}
\put(33.50,37.60){\oval(3,4.20)[b]}
\emline{35}{37.60}{7}{35}{38.65}{8}
\emline{31}{40.75}{9}{31}{38.65}{10}
\emline{34}{40.75}{11}{34}{39.70}{12}
\emline{33.33}{34.45}{13}{33.33}{35.50}{14}
\put(24,22.90){\circle{2}}
\emline{24}{26.05}{15}{24}{23.95}{16}
\emline{24}{21.85}{17}{24}{19.75}{18}
\emline{21}{26.05}{19}{21}{22.90}{20}
\put(21,23){\circle*{0.7}}
\put(43,23){\circle*{0.7}}
\put(31,8){\circle*{0.7}}
\put(33,8){\circle*{0.7}}
\put(12,8){\circle*{0.7}}
\put(3,26.05){\oval(4,6.30)[b]}
\emline{3}{19.75}{21}{3}{22.90}{22}
\put(44,22.90){\oval(2,4.20)[t]}
\put(43,22.90){\oval(4,4.20)[b]}
\emline{41}{26.05}{23}{41}{22.90}{24}
\emline{44}{26.05}{25}{44}{25}{26}
\emline{43}{19.75}{27}{43}{20.80}{28}
\put(34,8.20){\oval(2,4.20)[t]}
\emline{34}{11.35}{29}{34}{10.30}{30}
\emline{31}{11.35}{31}{31}{8.20}{32}
\emline{35}{5.05}{33}{35}{8.20}{34}
\emline{14}{5.05}{35}{14}{11.35}{36}
\emline{12}{11.35}{37}{12}{8.20}{38}
\put(5,28.15){\vector(-3,-2){0.2}}
\emline{11}{32.35}{39}{5}{28.15}{40}
\put(25,28.15){\vector(-3,-2){0.2}}
\emline{31}{32.35}{41}{25}{28.15}{42}
\put(41,28.15){\vector(3,-2){0.2}}
\emline{35}{32.35}{43}{41}{28.15}{44}
\emline{15}{32.35}{45}{16}{31.65}{46}
\emline{17}{30.95}{47}{18}{30.25}{48}
\put(12,13.45){\vector(3,-2){0.2}}
\emline{5}{17.65}{49}{12}{13.45}{50}
\put(14,13.45){\vector(-3,-2){0.2}}
\emline{21}{17.65}{51}{14}{13.45}{52}
\put(35,13.45){\vector(-3,-2){0.2}}
\emline{41}{17.65}{53}{35}{13.45}{54}
\put(7,30){\makebox(0,0)[rb]{$\id\times m$}}
\put(10,32){\makebox(0,0)[rb]{Leibniz}}
\put(17,30){\makebox(0,0)[rt]{$l$}}
\put(20,29){\makebox(0,0)[rt]{stochastic}}
\put(30,32){\makebox(0,0)[rb]{Leibniz}}
\put(27,30){\makebox(0,0)[rb]{$l\times\id$}}
\put(36,32){\makebox(0,0)[lb]{Borowiec}}
\put(39,30){\makebox(0,0)[lb]{$r\times\id$}}
\put(7.67,15.55){\makebox(0,0)[rt]{$l$}}
\put(16.67,15.55){\makebox(0,0)[rb]{$m$}}
\put(38.33,15.55){\makebox(0,0)[lt]{$l$}}
\put(19,37){\vector(-1,0){0.2}}
\emline{27}{37}{55}{19}{37}{56}
\put(19,8){\vector(-1,0){0.2}}
\emline{27}{8}{57}{19}{8}{58}
\put(23,9){\makebox(0,0)[cb]{$c$}}
\put(23,38){\makebox(0,0)[cb]{$a$}}
\put(31,13.44){\vector(4,-3){0.2}}
\emline{25}{18}{59}{31}{13.44}{60}
\put(28.22,16.22){\makebox(0,0)[lb]{$l$}}
\put(21,28){\vector(4,-3){0.2}}
\emline{19}{29.56}{61}{21}{28}{62}
\put(23,4){\makebox(0,0)[ct]{\begin{tabular}{c}\textsc{Figure} 11. The
Leibniz `heptagon' with 11 edges.\\An edge `$\id\times c$' is omitted
for simplicity.\end{tabular}}}\end{picture}$$}

\def\guM{\special{em:linewidth 0.4pt}\unitlength 1mm\linethickness{0.4pt}
$$\begin{picture}(62,26)          
\put(6,17){\oval(10,10)[b]}
\put(11.50,21){\oval(5,8)[b]}
\emline{1}{23}{1}{1}{17}{2}
\emline{5}{23}{3}{5}{21}{4}
\emline{5}{21}{5}{9}{21}{6}
\emline{14}{21}{7}{14}{23}{8}
\emline{7}{12}{9}{7}{10}{10}
\put(52,23.50){\oval(4,9)[b]}
\put(57,19){\oval(10,10)[b]}
\emline{62}{23}{11}{62}{19}{12}
\emline{57}{10}{13}{57}{14}{14}
\put(1,24){\makebox(0,0)[cb]{$L$}}
\put(6,9){\makebox(0,0)[ct]{$A$}}
\put(5,24){\makebox(0,0)[cb]{$A$}}
\put(14,24){\makebox(0,0)[cb]{$A$}}
\put(50,24){\makebox(0,0)[cb]{$L$}}
\put(54,24){\makebox(0,0)[cb]{$A$}}
\put(62,24){\makebox(0,0)[cb]{$A$}}
\put(15,21){\makebox(0,0)[lc]{$a$}}
\put(12,16){\makebox(0,0)[lt]{$m$}}
\put(7,13){\makebox(0,0)[cb]{$l$}}
\put(58,13){\makebox(0,0)[lt]{$m$}}
\put(51,18){\makebox(0,0)[rt]{$l$}}
\put(31,9){\makebox(0,0)[ct]{\begin{tabular}{c} The Leibniz operations.
\end{tabular}}}\end{picture}$$
$$\begin{picture}(100,20)   
\put(8,18){\oval(10,10)[b]}
\put(2.50,21){\oval(3,6)[b]}
\emline{13}{22}{1}{13}{18}{2}
\emline{8}{13}{3}{8}{10}{4}
\put(35.17,23){\oval(4.33,8)[b]}
\put(40,18.50){\oval(10,11)[b]}
\emline{45}{23}{5}{45}{19}{6}
\emline{40}{6}{7}{40}{13}{8}
\put(40,11){\circle*{1.33}}
\put(13,20){\circle*{1.33}}
\put(95,17){\circle{6}}
\put(95,14.50){\oval(4,5)[b]}
\emline{95}{9}{9}{95}{12}{10}
\emline{89}{23}{11}{89}{21}{12}
\emline{89}{21}{13}{93}{21}{14}
\emline{93}{21}{15}{93}{18.67}{16}
\emline{97}{23}{17}{97}{18.67}{18}
\put(92.17,23){\oval(14.33,12)[lb]}
\put(1.5,24){\makebox(0,0)[cb]{$L$}}
\put(4,24){\makebox(0,0)[cb]{$A$}}
\put(13,24){\makebox(0,0)[cb]{$A$}}
\put(3,15){\makebox(0,0)[rt]{$L$}}
\put(13,15){\makebox(0,0)[lt]{$A$}}
\put(8,9){\makebox(0,0)[ct]{$A$}}
\put(4.5,18){\makebox(0,0)[lc]{$r$}}
\put(15,20){\makebox(0,0)[lc]{$c$}}
\put(9,12){\makebox(0,0)[lt]{$l$}}
\put(33,19){\makebox(0,0)[rt]{$r$}}
\put(35,15){\makebox(0,0)[rt]{$L$}}
\put(40,14){\makebox(0,0)[cb]{$l$}}
\put(38,10){\makebox(0,0)[rc]{$c$}}
\put(95,17){\makebox(0,0)[cc]{$l$}}
\put(98,21){\makebox(0,0)[lc]{$a$}}
\put(96,11){\makebox(0,0)[lt]{$m$}}
\put(91,9){\makebox(0,0)[ct]{\begin{tabular}{c} A stochastic operation
\\is ternary.\end{tabular}}}
\put(22,6){\makebox(0,0)[ct]{\begin{tabular}{c} Borowiec's operations.
\end{tabular}}}
\put(50,1){\makebox(0,0)[ct]{\begin{tabular}{c} \textsc{Figure} 12.
\end{tabular}}}\end{picture}$$}

\def\guN{\special{em:linewidth 0.4pt}\unitlength 1mm\linethickness{0.4pt}
\begin{center}\begin{picture}(106,95)  
\put(49,10){\framebox(8,18)[cc]{}}
\put(73,10){\framebox(8,18)[cc]{}}
\emline{53}{26}{1}{53}{24}{2}
\emline{53}{22}{3}{53}{20}{4}
\put(53,16){\oval(4,4)[b]}
\put(53,18){\oval(4,4)[t]}
\emline{55}{18}{5}{55}{16}{6}
\emline{53}{12}{7}{53}{14}{8}
\put(53,24){\circle*{1}}
\put(53,22){\circle*{1}}
\put(51,18){\circle*{1}}
\put(51,16){\circle*{1}}
\put(1,38){\framebox(8,16)[cc]{}}
\put(25,38){\framebox(8,16)[cc]{}}
\put(49,38){\framebox(8,16)[cc]{}}
\put(73,38){\framebox(8,16)[cc]{}}
\put(97,38){\framebox(8,16)[cc]{}}
\put(49,64){\framebox(8,18)[cc]{}}
\put(73,64){\framebox(8,18)[cc]{}}
\put(53,70){\oval(4,4)[b]}
\put(53,72){\oval(4,4)[t]}
\emline{51}{72}{9}{51}{70}{10}
\emline{53}{66}{11}{53}{68}{12}
\emline{53}{74}{13}{53}{76}{14}
\put(53,76){\circle*{1}}
\put(55,72){\circle*{1}}
\put(55,70){\circle*{1}}
\emline{53}{80}{15}{53}{78}{16}
\put(53,78){\circle*{1}}
\put(77,71){\oval(4,6)[b]}
\put(77,75){\oval(4,6)[t]}
\emline{75}{75}{17}{75}{71}{18}
\emline{77}{66}{19}{77}{68}{20}
\emline{77}{80}{21}{77}{78}{22}
\put(77,17){\oval(4,6)[b]}
\put(77,21){\oval(4,6)[t]}
\emline{79}{21}{23}{79}{17}{24}
\emline{77}{26}{25}{77}{24}{26}
\emline{77}{12}{27}{77}{14}{28}
\put(79,75){\circle*{1.33}}
\put(79,71){\circle*{1.33}}
\put(75,17){\circle*{1.33}}
\put(75,21){\circle*{1.33}}
\emline{77}{52}{29}{77}{43}{30}
\put(77,40){\makebox(0,0)[cc]{sink}}
\put(101,46){\circle{4}}
\emline{101}{52}{31}{101}{48}{32}
\emline{101}{44}{33}{101}{40}{34}
\emline{53}{40}{35}{53}{44}{36}
\emline{53}{52}{37}{53}{48}{38}
\put(53,44){\circle*{1.33}}
\put(53,48){\circle*{1.33}}
\put(29,44){\circle{4}}
\emline{29}{40}{39}{29}{42}{40}
\emline{29}{46}{41}{29}{48}{42}
\emline{29}{52}{43}{29}{50}{44}
\put(29,48){\circle*{1}}
\put(29,50){\circle*{1}}
\put(5,44){\oval(4,4)[b]}
\put(5,48){\oval(4,4)[t]}
\emline{5}{52}{45}{5}{50}{46}
\emline{5}{40}{47}{5}{42}{48}
\put(3,44){\circle*{1.33}}
\put(7,44){\circle*{1.33}}
\put(7,48){\circle*{1.33}}
\put(3,48){\circle*{1.33}}
\put(29,37){\makebox(0,0)[ct]{source}}
\put(14,47){\vector(-1,-1){0.2}}
\emline{20}{47}{49}{18.50}{48.50}{50}
\emline{18.50}{48.50}{51}{15.50}{48.50}{52}
\emline{15.50}{48.50}{53}{14}{47}{54}
\put(20,45){\vector(1,1){0.2}}
\emline{14}{45}{55}{15.50}{43.50}{56}
\emline{15.50}{43.50}{57}{18.50}{43.50}{58}
\emline{18.50}{43.50}{59}{20}{45}{60}
\put(38,47){\vector(-1,-1){0.2}}
\emline{44}{47}{61}{42.50}{48.50}{62}
\emline{42.50}{48.50}{63}{39.50}{48.50}{64}
\emline{39.50}{48.50}{65}{38}{47}{66}
\put(62,47){\vector(-1,-1){0.2}}
\emline{68}{47}{67}{66.50}{48.50}{68}
\emline{66.50}{48.50}{69}{63.50}{48.50}{70}
\emline{63.50}{48.50}{71}{62}{47}{72}
\put(86,47){\vector(-1,-1){0.2}}
\emline{92}{47}{73}{90.50}{48.50}{74}
\emline{90.50}{48.50}{75}{87.50}{48.50}{76}
\emline{87.50}{48.50}{77}{86}{47}{78}
\put(44,45){\vector(1,1){0.2}}
\emline{38}{45}{79}{39.50}{43.50}{80}
\emline{39.50}{43.50}{81}{42.50}{43.50}{82}
\emline{42.50}{43.50}{83}{44}{45}{84}
\put(68,45){\vector(1,1){0.2}}
\emline{62}{45}{85}{63.50}{43.50}{86}
\emline{63.50}{43.50}{87}{66.50}{43.50}{88}
\emline{66.50}{43.50}{89}{68}{45}{90}
\put(92,45){\vector(1,1){0.2}}
\emline{86}{45}{91}{87.50}{43.50}{92}
\emline{87.50}{43.50}{93}{90.50}{43.50}{94}
\emline{90.50}{43.50}{95}{92}{45}{96}
\put(62,74){\vector(-1,-1){0.2}}
\emline{68}{74}{97}{66.50}{75.50}{98}
\emline{66.50}{75.50}{99}{63.50}{75.50}{100}
\emline{63.50}{75.50}{101}{62}{74}{102}
\put(68,72){\vector(1,1){0.2}}
\emline{62}{72}{103}{63.50}{70.50}{104}
\emline{63.50}{70.50}{105}{66.50}{70.50}{106}
\emline{66.50}{70.50}{107}{68}{72}{108}
\put(62,20){\vector(-1,-1){0.2}}
\emline{68}{20}{109}{66.50}{21.50}{110}
\emline{66.50}{21.50}{111}{63.50}{21.50}{112}
\emline{63.50}{21.50}{113}{62}{20}{114}
\put(68,18){\vector(1,1){0.2}}
\emline{62}{18}{115}{63.50}{16.50}{116}
\emline{63.50}{16.50}{117}{66.50}{16.50}{118}
\emline{66.50}{16.50}{119}{68}{18}{120}
\put(52,30){\vector(1,-3){0.2}}
\emline{52}{36}{121}{51.08}{33.86}{122}
\emline{51.08}{33.86}{123}{52}{30}{124}
\put(54,36){\vector(-1,3){0.2}}
\emline{54}{30}{125}{54.92}{32.14}{126}
\emline{54.92}{32.14}{127}{54}{36}{128}
\put(76,30){\vector(1,-3){0.2}}
\emline{76}{36}{129}{75.08}{33.86}{130}
\emline{75.08}{33.86}{131}{76}{30}{132}
\put(78,36){\vector(-1,3){0.2}}
\emline{78}{30}{133}{78.92}{32.14}{134}
\emline{78.92}{32.14}{135}{78}{36}{136}
\put(52,56){\vector(1,-3){0.2}}
\emline{52}{62}{137}{51.08}{59.86}{138}
\emline{51.08}{59.86}{139}{52}{56}{140}
\put(54,62){\vector(-1,3){0.2}}
\emline{54}{56}{141}{54.92}{58.14}{142}
\emline{54.92}{58.14}{143}{54}{62}{144}
\put(76,56){\vector(1,-3){0.2}}
\emline{76}{62}{145}{75.08}{59.86}{146}
\emline{75.08}{59.86}{147}{76}{56}{148}
\put(78,62){\vector(-1,3){0.2}}
\emline{78}{56}{149}{78.92}{58.14}{150}
\emline{78.92}{58.14}{151}{78}{62}{152}
\emline{5}{33}{153}{7.77}{30.62}{154}
\emline{7.77}{30.62}{155}{10.54}{28.36}{156}
\emline{10.54}{28.36}{157}{13.31}{26.20}{158}
\emline{13.31}{26.20}{159}{16.08}{24.16}{160}
\emline{16.08}{24.16}{161}{18.85}{22.24}{162}
\emline{18.85}{22.24}{163}{21.62}{20.43}{164}
\emline{21.62}{20.43}{165}{24.38}{18.73}{166}
\emline{24.38}{18.73}{167}{27.15}{17.14}{168}
\emline{27.15}{17.14}{169}{29.92}{15.67}{170}
\emline{29.92}{15.67}{171}{32.69}{14.31}{172}
\emline{32.69}{14.31}{173}{35.46}{13.07}{174}
\emline{35.46}{13.07}{175}{38.23}{11.94}{176}
\emline{38.23}{11.94}{177}{41}{10.92}{178}
\emline{41}{10.92}{179}{43.77}{10.01}{180}
\emline{43.77}{10.01}{181}{46.54}{9.22}{182}
\emline{46.54}{9.22}{183}{49.31}{8.54}{184}
\emline{49.31}{8.54}{185}{52.08}{7.98}{186}
\emline{52.08}{7.98}{187}{54.85}{7.53}{188}
\emline{54.85}{7.53}{189}{57.61}{7.19}{190}
\emline{57.61}{7.19}{191}{60.38}{6.96}{192}
\emline{60.38}{6.96}{193}{63.15}{6.85}{194}
\emline{63.15}{6.85}{195}{65.92}{6.86}{196}
\emline{65.92}{6.86}{197}{68.69}{6.97}{198}
\emline{68.69}{6.97}{199}{71.46}{7.20}{200}
\emline{71.46}{7.20}{201}{74.23}{7.54}{202}
\emline{74.23}{7.54}{203}{77}{8}{204}
\emline{5}{59}{205}{7.77}{61.40}{206}
\emline{7.77}{61.40}{207}{10.54}{63.69}{208}
\emline{10.54}{63.69}{209}{13.31}{65.87}{210}
\emline{13.31}{65.87}{211}{16.08}{67.92}{212}
\emline{16.08}{67.92}{213}{18.85}{69.87}{214}
\emline{18.85}{69.87}{215}{21.62}{71.69}{216}
\emline{21.62}{71.69}{217}{24.38}{73.40}{218}
\emline{24.38}{73.40}{219}{27.15}{75}{220}
\emline{27.15}{75}{221}{29.92}{76.48}{222}
\emline{29.92}{76.48}{223}{32.69}{77.85}{224}
\emline{32.69}{77.85}{225}{35.46}{79.10}{226}
\emline{35.46}{79.10}{227}{38.23}{80.23}{228}
\emline{38.23}{80.23}{229}{41}{81.25}{230}
\emline{41}{81.25}{231}{43.77}{82.15}{232}
\emline{43.77}{82.15}{233}{46.54}{82.94}{234}
\emline{46.54}{82.94}{235}{49.31}{83.62}{236}
\emline{49.31}{83.62}{237}{52.08}{84.17}{238}
\emline{52.08}{84.17}{239}{54.85}{84.62}{240}
\emline{54.85}{84.62}{241}{57.61}{84.94}{242}
\emline{57.61}{84.94}{243}{60.38}{85.15}{244}
\emline{60.38}{85.15}{245}{63.15}{85.25}{246}
\emline{63.15}{85.25}{247}{65.92}{85.23}{248}
\emline{65.92}{85.23}{249}{68.69}{85.10}{250}
\emline{68.69}{85.10}{251}{71.46}{84.85}{252}
\emline{71.46}{84.85}{253}{74.23}{84.48}{254}
\emline{74.23}{84.48}{255}{77}{84}{256}
\emline{29}{59}{257}{31.13}{62.56}{258}
\emline{31.13}{62.56}{259}{33.26}{65.91}{260}
\emline{33.26}{65.91}{261}{35.40}{69.04}{262}
\emline{35.40}{69.04}{263}{37.53}{71.95}{264}
\emline{37.53}{71.95}{265}{39.66}{74.64}{266}
\emline{39.66}{74.64}{267}{41.79}{77.12}{268}
\emline{41.79}{77.12}{269}{43.93}{79.38}{270}
\emline{43.93}{79.38}{271}{46.06}{81.42}{272}
\emline{46.06}{81.42}{273}{48.19}{83.24}{274}
\emline{48.19}{83.24}{275}{50.32}{84.85}{276}
\emline{50.32}{84.85}{277}{52.46}{86.24}{278}
\emline{52.46}{86.24}{279}{54.59}{87.41}{280}
\emline{54.59}{87.41}{281}{56.72}{88.36}{282}
\emline{56.72}{88.36}{283}{58.85}{89.10}{284}
\emline{58.85}{89.10}{285}{60.99}{89.61}{286}
\emline{60.99}{89.61}{287}{63.12}{89.92}{288}
\emline{63.12}{89.92}{289}{65.25}{90}{290}
\emline{65.25}{90}{291}{67.38}{89.86}{292}
\emline{67.38}{89.86}{293}{69.52}{89.51}{294}
\emline{69.52}{89.51}{295}{71.65}{88.94}{296}
\emline{71.65}{88.94}{297}{73.78}{88.16}{298}
\emline{73.78}{88.16}{299}{75.91}{87.15}{300}
\emline{75.91}{87.15}{301}{78.05}{85.93}{302}
\emline{78.05}{85.93}{303}{80.18}{84.49}{304}
\emline{80.18}{84.49}{305}{82.31}{82.83}{306}
\emline{82.31}{82.83}{307}{84.44}{80.96}{308}
\emline{84.44}{80.96}{309}{86.58}{78.86}{310}
\emline{86.58}{78.86}{311}{88.71}{76.55}{312}
\emline{88.71}{76.55}{313}{90.84}{74.03}{314}
\emline{90.84}{74.03}{315}{92.97}{71.28}{316}
\emline{92.97}{71.28}{317}{95.11}{68.32}{318}
\emline{95.11}{68.32}{319}{97.24}{65.14}{320}
\emline{97.24}{65.14}{321}{101}{59}{322}
\put(65,77){\makebox(0,0)[cb]{$\delta$}}
\put(65,69){\makebox(0,0)[ct]{$d$}}
\put(65,50){\makebox(0,0)[cb]{$\delta$}}
\put(65,42){\makebox(0,0)[ct]{$d$}}
\put(65,23){\makebox(0,0)[cb]{$\delta$}}
\put(65,15){\makebox(0,0)[ct]{$d$}}
\put(41,50){\makebox(0,0)[cb]{$\triangle$}}
\put(41,42){\makebox(0,0)[ct]{$m$}}
\put(89,42){\makebox(0,0)[ct]{$\triangle$}}
\put(89,50){\makebox(0,0)[cb]{$m$}}
\put(50,59){\makebox(0,0)[rc]{$m_r$}}
\put(56,59){\makebox(0,0)[lc]{$\triangle_r$}}
\put(74,59){\makebox(0,0)[rc]{$m_r$}}
\put(80,59){\makebox(0,0)[lc]{$\triangle_r$}}
\put(50,33){\makebox(0,0)[rc]{$\triangle_l$}}
\put(56,33){\makebox(0,0)[lc]{$m_l$}}
\put(74,33){\makebox(0,0)[rc]{$\triangle_l$}}
\put(80,33){\makebox(0,0)[lc]{$m_l$}}
\put(88,80){\makebox(0,0)[lc]{stochastic \&}}
\put(91,77){\makebox(0,0)[lc]{costochastic}}
\put(94,73){\makebox(0,0)[lc]{(co)differential}}
\put(98,69){\makebox(0,0)[lc]{$d$ \& $\delta$}}
\put(15,68){\makebox(0,0)[rc]{$d\times\id\nearrow$}}
\put(15,66){\makebox(0,0)[lc]{$\swarrow\delta\times\id$}}
\put(16,26){\makebox(0,0)[lc]{$\searrow\id\times d$}}
\put(15,23){\makebox(0,0)[rc]{$\id\times\delta\nwarrow$}}
\put(53,5){\makebox(0,0)[ct]{\begin{tabular}{c}\textsc{Figure} 13.
A fragment of the Cartan garden $1\mapsto 1$ with plants in the
shelters.\end{tabular}}}\end{picture}\end{center}}

\def\guO{\special{em:linewidth 0.4pt}\unitlength 1mm\linethickness{0.4pt}
\begin{center}\begin{picture}(92,34)       
\put(4,24){\circle{4}}
\put(88,24){\circle{4}}
\put(90,24){\circle*{1.3}}
\put(2,24){\circle*{1.3}}
\put(46,22.50){\circle*{1.3}}
\put(46,25.50){\circle*{1.3}}
\emline{4}{28}{1}{4}{26}{2}
\emline{4}{22}{3}{4}{20}{4}
\emline{46}{20}{5}{46}{23}{6}
\emline{46}{28}{7}{46}{25}{8}
\emline{88}{28}{9}{88}{26}{10}
\emline{88}{20}{11}{88}{22}{12}
\put(41,25){\vector(4,-1){0.2}}
\emline{11}{25}{13}{13.42}{25.59}{14}
\emline{13.42}{25.59}{15}{15.84}{26.08}{16}
\emline{15.84}{26.08}{17}{18.26}{26.47}{18}
\emline{18.26}{26.47}{19}{20.68}{26.75}{20}
\emline{20.68}{26.75}{21}{23.10}{26.93}{22}
\emline{23.10}{26.93}{23}{25.52}{27}{24}
\emline{25.52}{27}{25}{27.94}{26.97}{26}
\emline{27.94}{26.97}{27}{30.35}{26.83}{28}
\emline{30.35}{26.83}{29}{32.77}{26.59}{30}
\emline{32.77}{26.59}{31}{35.19}{26.25}{32}
\emline{35.19}{26.25}{33}{37.61}{25.80}{34}
\emline{37.61}{25.80}{35}{41}{25}{36}
\put(11,23){\vector(-4,1){0.2}}
\emline{41}{23}{37}{38.58}{22.41}{38}
\emline{38.58}{22.41}{39}{36.16}{21.92}{40}
\emline{36.16}{21.92}{41}{33.74}{21.53}{42}
\emline{33.74}{21.53}{43}{31.32}{21.25}{44}
\emline{31.32}{21.25}{45}{28.90}{21.07}{46}
\emline{28.90}{21.07}{47}{26.48}{21}{48}
\emline{26.48}{21}{49}{24.06}{21.03}{50}
\emline{24.06}{21.03}{51}{21.65}{21.17}{52}
\emline{21.65}{21.17}{53}{19.23}{21.41}{54}
\emline{19.23}{21.41}{55}{16.81}{21.75}{56}
\emline{16.81}{21.75}{57}{14.39}{22.20}{58}
\emline{14.39}{22.20}{59}{11}{23}{60}
\put(81,25){\vector(4,-1){0.2}}
\emline{51}{25}{61}{53.42}{25.59}{62}
\emline{53.42}{25.59}{63}{55.84}{26.08}{64}
\emline{55.84}{26.08}{65}{58.26}{26.47}{66}
\emline{58.26}{26.47}{67}{60.68}{26.75}{68}
\emline{60.68}{26.75}{69}{63.10}{26.93}{70}
\emline{63.10}{26.93}{71}{65.52}{27}{72}
\emline{65.52}{27}{73}{67.94}{26.97}{74}
\emline{67.94}{26.97}{75}{70.35}{26.83}{76}
\emline{70.35}{26.83}{77}{72.77}{26.59}{78}
\emline{72.77}{26.59}{79}{75.19}{26.25}{80}
\emline{75.19}{26.25}{81}{77.61}{25.80}{82}
\emline{77.61}{25.80}{83}{81}{25}{84}
\put(51,23){\vector(-4,1){0.2}}
\emline{81}{23}{85}{78.58}{22.41}{86}
\emline{78.58}{22.41}{87}{76.16}{21.92}{88}
\emline{76.16}{21.92}{89}{73.74}{21.53}{90}
\emline{73.74}{21.53}{91}{71.32}{21.25}{92}
\emline{71.32}{21.25}{93}{68.90}{21.07}{94}
\emline{68.90}{21.07}{95}{66.48}{21}{96}
\emline{66.48}{21}{97}{64.06}{21.03}{98}
\emline{64.06}{21.03}{99}{61.65}{21.17}{100}
\emline{61.65}{21.17}{101}{59.23}{21.41}{102}
\emline{59.23}{21.41}{103}{56.81}{21.75}{104}
\emline{56.81}{21.75}{105}{54.39}{22.20}{106}
\emline{54.39}{22.20}{107}{51}{23}{108}
\put(26,28){\makebox(0,0)[cb]{left evaluation}}
\put(66,28){\makebox(0,0)[cb]{left coevaluation}}
\put(26,20){\makebox(0,0)[ct]{right coevaluation}}
\put(66,20){\makebox(0,0)[ct]{right evaluation}}
\put(46,15){\makebox(0,0)[ct]
{\begin{tabular}{c}\textsc{Figure} 14. Two distinct closed structures:
a plant $\guj$ is said to be\\a left or a right adjoint wrt to unit
creator $\in(0\mapsto 1).$\end{tabular}}}\end{picture}\end{center}}

\def\guP{\special{em:linewidth 0.4pt}\unitlength 1mm\linethickness{0.4pt}
\begin{center}\begin{picture}(120,52)       
\put(4,22){\oval(4,4)[b]}
\put(4,26){\oval(4,4)[t]}
\emline{6}{22}{1}{6}{26}{2}
\emline{4}{30}{3}{4}{28}{4}
\emline{4}{18}{5}{4}{20}{6}
\put(2,22){\circle*{1.33}}
\put(2,26){\circle*{1.33}}
\emline{23}{8}{7}{23}{16}{8}
\emline{97}{8}{9}{97}{16}{10}
\put(97,13){\circle*{1.33}}
\put(42,22){\oval(4,4)[b]}
\put(42,26){\oval(4,4)[t]}
\emline{40}{22}{11}{40}{26}{12}
\emline{42}{18}{13}{42}{20}{14}
\emline{42}{30}{15}{42}{28}{16}
\put(44,26){\circle*{1.33}}
\put(44,22){\circle*{1.33}}
\put(11,46){\circle{4}}
\put(35,46){\circle{4}}
\put(13,48){\oval(4,4)[t]}
\put(13,44){\oval(4,4)[b]}
\emline{15}{44}{17}{15}{48}{18}
\emline{13}{50}{19}{13}{52}{20}
\emline{13}{40}{21}{13}{42}{22}
\put(33,44){\oval(4,4)[b]}
\put(33,48){\oval(4,4)[t]}
\emline{31}{44}{23}{31}{48}{24}
\emline{33}{52}{25}{33}{50}{26}
\emline{33}{40}{27}{33}{42}{28}
\put(78,22){\oval(4,4)[b]}
\put(78,26){\oval(4,4)[t]}
\emline{76}{22}{29}{76}{26}{30}
\emline{78}{18}{31}{78}{20}{32}
\emline{78}{30}{33}{78}{28}{34}
\put(80,26){\circle*{1.33}}
\put(80,22){\circle*{1.33}}
\put(116,22){\oval(4,4)[b]}
\put(116,26){\oval(4,4)[t]}
\emline{118}{22}{35}{118}{26}{36}
\emline{116}{30}{37}{116}{28}{38}
\emline{116}{18}{39}{116}{20}{40}
\put(114,22){\circle*{1.33}}
\put(114,26){\circle*{1.33}}
\put(107,46){\circle{4}}
\put(109,48){\oval(4,4)[t]}
\put(109,44){\oval(4,4)[b]}
\emline{111}{44}{41}{111}{48}{42}
\emline{109}{50}{43}{109}{52}{44}
\emline{109}{40}{45}{109}{42}{46}
\put(87,46){\circle{4}}
\put(85,44){\oval(4,4)[b]}
\put(85,48){\oval(4,4)[t]}
\emline{83}{44}{47}{83}{48}{48}
\emline{85}{52}{49}{85}{50}{50}
\emline{85}{40}{51}{85}{42}{52}
\put(13,46){\circle*{1.33}}
\put(33,46){\circle*{1.33}}
\put(83,46){\circle*{1.33}}
\put(89,46){\circle*{1.33}}
\put(105,46){\circle*{1.33}}
\put(111,46){\circle*{1.33}}
\put(76,24){\circle*{1.33}}
\put(118,24){\circle*{1.33}}
\put(27,46){\vector(1,0){0.2}}
\emline{19}{46}{53}{27}{46}{54}
\put(101,46){\vector(1,0){0.2}}
\emline{93}{46}{55}{101}{46}{56}
\put(9,39){\vector(1,2){0.2}}
\emline{5}{32}{57}{9}{39}{58}
\put(83,39){\vector(1,2){0.2}}
\emline{79}{32}{59}{83}{39}{60}
\put(8.33,36){\makebox(0,0)[lt]{$\coev_l$}}
\put(82,36){\makebox(0,0)[lt]{$\coev_l$}}
\put(23,47){\makebox(0,0)[cb]{$a$}}
\put(97,47){\makebox(0,0)[cb]{$a^*$}}
\put(41,32){\vector(1,-2){0.2}}
\emline{37}{39}{61}{41}{32}{62}
\put(115,32){\vector(1,-2){0.2}}
\emline{111}{39}{63}{115}{32}{64}
\put(40,36){\makebox(0,0)[lb]{$\ev_l$}}
\put(114,36){\makebox(0,0)[lb]{$\ev_l$}}
\put(29,14){\vector(-1,-1){0.2}}
\emline{35}{20}{65}{29}{14}{66}
\put(103,14){\vector(-1,-1){0.2}}
\emline{109}{20}{67}{103}{14}{68}
\put(85,20){\vector(-1,1){0.2}}
\emline{91}{14}{69}{85}{20}{70}
\put(11,20){\vector(-1,1){0.2}}
\emline{17}{14}{71}{11}{20}{72}
\put(67,24){\vector(1,0){0.2}}
\emline{53}{24}{73}{67}{24}{74}
\put(60,25){\makebox(0,0)[cb]{geometry}}
\put(57,5){\makebox(0,0)[ct]{\begin{tabular}{c}\textsc{Figure} 15. The
pentagons of evaluation for a reflexive left adjoint.\end{tabular}}}
\end{picture}\end{center}}

\def\guR{\special{em:linewidth 0.4pt}\unitlength 1mm\linethickness{0.4pt}
\begin{center}\begin{picture}(74,15)       
\emline{2}{5}{1}{2}{8}{2}
\emline{2}{10}{3}{2}{13}{4}
\put(2,10.67){\circle*{1.2}}
\put(2,7.67){\circle*{1.2}}
\emline{19}{5}{5}{19}{13}{6}
\emline{55}{5}{7}{55}{13}{8}
\put(55,9){\circle*{1.3}}
\emline{72}{5}{9}{72}{7}{10}
\emline{72}{10.5}{11}{72}{13}{12}
\put(72,10.5){\circle*{1.2}}
\put(72,7.67){\makebox(0,0)[cc]{$\star$}}
\put(68,9.67){\vector(2,-1){0.2}}
\emline{59}{9.67}{13}{61.46}{10.60}{14}
\emline{61.46}{10.60}{15}{63.67}{10.99}{16}
\emline{63.67}{10.99}{17}{65.64}{10.82}{18}
\emline{65.64}{10.82}{19}{68}{9.67}{20}
\put(59,8.33){\vector(-2,1){0.2}}
\emline{68}{8.33}{21}{65.44}{7.33}{22}
\emline{65.44}{7.33}{23}{60.94}{7.33}{24}
\emline{60.94}{7.33}{25}{59}{8.33}{26}
\put(15,9.67){\vector(2,-1){0.2}}
\emline{6}{9.67}{27}{8.46}{10.60}{28}
\emline{8.46}{10.60}{29}{10.67}{10.99}{30}
\emline{10.67}{10.99}{31}{12.64}{10.82}{32}
\emline{12.64}{10.82}{33}{15}{9.67}{34}
\put(6,8.33){\vector(-2,1){0.2}}
\emline{15}{8.33}{35}{12.44}{7.33}{36}
\emline{12.44}{7.33}{37}{7.94}{7.33}{38}
\emline{7.94}{7.33}{39}{6}{8.33}{40}
\put(10,12){\makebox(0,0)[cb]{$d$}}
\put(11,6){\makebox(0,0)[ct]{$\delta$}}
\put(63,12){\makebox(0,0)[cb]{$d^*$}}
\put(64,6){\makebox(0,0)[ct]{$\delta^*$}}
\put(44,9){\vector(1,0){0.2}}
\emline{30}{9}{41}{44}{9}{42}
\put(37,10){\makebox(0,0)[cb]{geometry}}
\put(37,3){\makebox(0,0)[ct]{\begin{tabular}{c}\textsc{Figure} 16.
The Cartan paths on the left and the translators.\end{tabular}}}
\end{picture}\end{center}}

\def\guS{\special{em:linewidth 0.4pt}\unitlength 1mm\linethickness{0.4pt}
\begin{center}\begin{picture}(78,41)       
\put(1,5){\framebox(8,16)[cc]{}}
\put(37,5){\framebox(4,16)[cc]{}}
\put(69,5){\framebox(8,16)[cc]{}}
\put(23,35){\circle{4}}
\put(55,35){\circle{4}}
\emline{55}{39}{1}{55}{37}{2}
\emline{55}{33}{3}{55}{31}{4}
\emline{23}{31}{5}{23}{33}{6}
\emline{23}{39}{7}{23}{37}{8}
\put(21,35){\circle*{1.33}}
\put(57,35){\circle*{1.33}}
\emline{39}{19}{9}{39}{15}{10}
\put(39,15){\circle*{1.33}}
\emline{39}{7}{11}{39}{11}{12}
\put(39,11){\circle*{1.33}}
\put(5,11){\oval(4,4)[b]}
\put(5,15){\oval(4,4)[t]}
\put(73,11){\oval(4,4)[b]}
\put(73,15){\oval(4,4)[t]}
\emline{73}{19}{13}{73}{17}{14}
\emline{73}{7}{15}{73}{9}{16}
\emline{5}{7}{17}{5}{9}{18}
\emline{5}{19}{19}{5}{17}{20}
\put(3,15){\circle*{1.33}}
\put(7,15){\circle*{1.33}}
\put(71,15){\circle*{1.33}}
\put(75,15){\circle*{1.33}}
\put(71,11){\circle*{1.33}}
\put(7,11){\circle*{1.33}}
\put(75,12){\makebox(0,0)[cc]{$\star$}}
\put(3,12){\makebox(0,0)[cc]{$\star$}}
\put(28,14){\vector(2,-1){0.2}}
\emline{18}{14}{21}{19.92}{15.24}{22}
\emline{19.92}{15.24}{23}{21.85}{15.89}{24}
\emline{21.85}{15.89}{25}{23.77}{15.95}{26}
\emline{23.77}{15.95}{27}{25.69}{15.42}{28}
\emline{25.69}{15.42}{29}{28}{14}{30}
\put(18,12){\vector(-2,1){0.2}}
\emline{28}{12}{31}{26.08}{10.76}{32}
\emline{26.08}{10.76}{33}{24.15}{10.11}{34}
\emline{24.15}{10.11}{35}{22.23}{10.05}{36}
\emline{22.23}{10.05}{37}{20.31}{10.58}{38}
\emline{20.31}{10.58}{39}{18}{12}{40}
\put(50,14){\vector(-2,-1){0.2}}
\emline{60}{14}{41}{58.08}{15.24}{42}
\emline{58.08}{15.24}{43}{56.15}{15.89}{44}
\emline{56.15}{15.89}{45}{54.23}{15.95}{46}
\emline{54.23}{15.95}{47}{52.31}{15.42}{48}
\emline{52.31}{15.42}{49}{50}{14}{50}
\put(60,12){\vector(2,1){0.2}}
\emline{50}{12}{51}{51.92}{10.76}{52}
\emline{51.92}{10.76}{53}{53.85}{10.11}{54}
\emline{53.85}{10.11}{55}{55.77}{10.05}{56}
\emline{55.77}{10.05}{57}{57.69}{10.58}{58}
\emline{57.69}{10.58}{59}{60}{12}{60}
\put(23,17){\makebox(0,0)[cb]{$l$}}
\put(55,17){\makebox(0,0)[cb]{$r$}}
\put(55,9){\makebox(0,0)[ct]{$R$}}
\put(23,9){\makebox(0,0)[ct]{$L$}}
\put(16,35){\vector(1,0){0.2}}
\emline{5}{26}{61}{5.32}{28.18}{62}
\emline{5.32}{28.18}{63}{5.94}{30.06}{64}
\emline{5.94}{30.06}{65}{6.86}{31.64}{66}
\emline{6.86}{31.64}{67}{8.08}{32.92}{68}
\emline{8.08}{32.92}{69}{9.61}{33.89}{70}
\emline{9.61}{33.89}{71}{11.44}{34.56}{72}
\emline{11.44}{34.56}{73}{13.57}{34.93}{74}
\emline{13.57}{34.93}{75}{16}{35}{76}
\put(62,35){\vector(-1,0){0.2}}
\emline{73}{26}{77}{72.58}{28.29}{78}
\emline{72.58}{28.29}{79}{71.85}{30.24}{80}
\emline{71.85}{30.24}{81}{70.81}{31.86}{82}
\emline{70.81}{31.86}{83}{69.45}{33.15}{84}
\emline{69.45}{33.15}{85}{67.79}{34.10}{86}
\emline{67.79}{34.10}{87}{65.81}{34.71}{88}
\emline{65.81}{34.71}{89}{62}{35}{90}
\put(38,26){\vector(1,-4){0.2}}
\emline{28}{35}{91}{30.38}{34.92}{92}
\emline{30.38}{34.92}{93}{32.44}{34.51}{94}
\emline{32.44}{34.51}{95}{34.18}{33.76}{96}
\emline{34.18}{33.76}{97}{35.59}{32.67}{98}
\emline{35.59}{32.67}{99}{36.68}{31.25}{100}
\emline{36.68}{31.25}{101}{37.45}{29.50}{102}
\emline{37.45}{29.50}{103}{38}{26}{104}
\put(40,26){\vector(-1,-4){0.2}}
\emline{50}{35}{105}{47.62}{34.92}{106}
\emline{47.62}{34.92}{107}{45.56}{34.51}{108}
\emline{45.56}{34.51}{109}{43.82}{33.76}{110}
\emline{43.82}{33.76}{111}{42.41}{32.67}{112}
\emline{42.41}{32.67}{113}{41.32}{31.25}{114}
\emline{41.32}{31.25}{115}{40.55}{29.50}{116}
\emline{40.55}{29.50}{117}{40}{26}{118}
\put(7,25){\vector(0,-1){0.2}}
\emline{17}{31}{119}{14.47}{31.43}{120}
\emline{14.47}{31.43}{121}{12.30}{31.49}{122}
\emline{12.30}{31.49}{123}{10.50}{31.18}{124}
\emline{10.50}{31.18}{125}{9.08}{30.50}{126}
\emline{9.08}{30.50}{127}{8.02}{29.45}{128}
\emline{8.02}{29.45}{129}{7.33}{28.04}{130}
\emline{7.33}{28.04}{131}{7}{25}{132}
\put(27,33){\vector(-1,0){0.2}}
\emline{36}{25}{133}{35.60}{27.13}{134}
\emline{35.60}{27.13}{135}{34.83}{28.93}{136}
\emline{34.83}{28.93}{137}{33.69}{30.41}{138}
\emline{33.69}{30.41}{139}{32.17}{31.56}{140}
\emline{32.17}{31.56}{141}{30.28}{32.39}{142}
\emline{30.28}{32.39}{143}{27}{33}{144}
\put(51,33){\vector(4,1){0.2}}
\emline{42}{25}{145}{42.53}{27.26}{146}
\emline{42.53}{27.26}{147}{43.44}{29.15}{148}
\emline{43.44}{29.15}{149}{44.75}{30.67}{150}
\emline{44.75}{30.67}{151}{46.44}{31.81}{152}
\emline{46.44}{31.81}{153}{48.53}{32.59}{154}
\emline{48.53}{32.59}{155}{51}{33}{156}
\put(71,25){\vector(0,-1){0.2}}
\emline{61}{31}{157}{63.53}{31.43}{158}
\emline{63.53}{31.43}{159}{65.70}{31.49}{160}
\emline{65.70}{31.49}{161}{67.50}{31.18}{162}
\emline{67.50}{31.18}{163}{68.92}{30.50}{164}
\emline{68.92}{30.50}{165}{69.98}{29.45}{166}
\emline{69.98}{29.45}{167}{70.67}{28.04}{168}
\emline{70.67}{28.04}{169}{71}{25}{170}
\put(7,34){\makebox(0,0)[rc]{$\delta^*\times d$}}
\put(71,34){\makebox(0,0)[lc]{$d\times\delta^*$}}
\put(69,30){\makebox(0,0)[rt]{$\delta\times d^*$}}
\put(10,30){\makebox(0,0)[lt]{$d^*\times\delta$}}
\put(32,36){\makebox(0,0)[cb]{ev$_l$}}
\put(46,36){\makebox(0,0)[cb]{ev$_r$}}
\put(33,30){\makebox(0,0)[rt]{coev$_r$}}
\put(46,30){\makebox(0,0)[lt]{coev$_l$}}
\put(39,3){\makebox(0,0)[ct]{\begin{tabular}{c}\textsc{Figure} 17.
The Leibniz 2-cells from Cartan.\end{tabular}}}\end{picture}\end{center}}

\def\guT{\special{em:linewidth 0.4pt}\unitlength 1mm\linethickness{0.4pt}
\begin{center}\begin{picture}(92,51)        
\put(57,5){\framebox(4,16)[cc]{}}
\put(83,5){\framebox(8,16)[cc]{}}
\put(10,31){\framebox(10,20)[cc]{}}
\put(41,31){\framebox(10,20)[cc]{}}
\put(72,31){\framebox(10,20)[cc]{}}
\put(14,41){\circle{4}}
\put(16,41){\circle*{1.3}}
\put(1,5){\framebox(8,16)[cc]{}}
\put(31,5){\framebox(4,16)[cc]{}}
\put(16,39){\oval(4,8)[b]}
\put(16,43){\oval(4,8)[t]}
\put(18,43){\circle*{1.3}}
\put(18,39){\circle*{1.3}}
\emline{16}{49}{1}{16}{47}{2}
\emline{16}{35}{3}{16}{33}{4}
\put(45,37){\oval(4,4)[b]}
\put(47,39){\oval(4,4)[b]}
\put(45,45){\oval(4,4)[t]}
\put(47,43){\oval(4,4)[t]}
\emline{43}{45}{5}{43}{37}{6}
\emline{45}{43}{7}{45}{39}{8}
\emline{45}{49}{9}{45}{47}{10}
\emline{45}{35}{11}{45}{33}{12}
\put(45,41){\circle*{1.3}}
\put(49,43){\circle*{1.3}}
\put(49,39){\circle*{1.3}}
\put(76,37){\oval(4,4)[b]}
\put(78,39){\oval(4,4)[b]}
\put(76,45){\oval(4,4)[t]}
\put(78,43){\oval(4,4)[t]}
\emline{74}{45}{13}{74}{37}{14}
\emline{76}{49}{15}{76}{47}{16}
\emline{76}{33}{17}{76}{35}{18}
\put(76,43){\circle*{1.3}}
\put(80,43){\circle*{1.3}}
\put(80,39){\circle*{1.3}}
\put(76,39){\makebox(0,0)[cc]{$\star$}}
\put(5,11){\oval(4,4)[b]}
\put(5,15){\oval(4,4)[t]}
\emline{5}{19}{19}{5}{17}{20}
\emline{5}{7}{21}{5}{9}{22}
\put(3,11){\circle*{1.3}}
\put(7,11){\circle*{1.3}}
\put(7,15){\circle*{1.3}}
\put(3,15){\circle*{1.3}}
\put(87,11){\oval(4,4)[b]}
\put(87,15){\oval(4,4)[t]}
\emline{85}{15}{23}{85}{11}{24}
\emline{87}{19}{25}{87}{17}{26}
\emline{87}{7}{27}{87}{9}{28}
\put(89,15){\circle*{1.3}}
\put(89,11){\circle*{1.3}}
\emline{33}{19}{29}{33}{15}{30}
\put(33,15){\circle*{1.3}}
\emline{33}{7}{31}{33}{11}{32}
\put(33,11){\circle*{1.3}}
\emline{59}{7}{33}{59}{19}{34}
\put(34,42){\vector(4,-3){0.2}}
\emline{27}{42}{35}{28.53}{43.40}{36}
\emline{28.53}{43.40}{37}{30.10}{43.98}{38}
\emline{30.10}{43.98}{39}{31.70}{43.74}{40}
\emline{31.70}{43.74}{41}{34}{42}{42}
\put(27,40){\vector(-4,3){0.2}}
\emline{34}{40}{43}{32.35}{38.60}{44}
\emline{32.35}{38.60}{45}{30.74}{38.02}{46}
\emline{30.74}{38.02}{47}{29.15}{38.26}{48}
\emline{29.15}{38.26}{49}{27}{40}{50}
\put(65,42){\vector(4,-3){0.2}}
\emline{58}{42}{51}{59.53}{43.40}{52}
\emline{59.53}{43.40}{53}{61.10}{43.98}{54}
\emline{61.10}{43.98}{55}{62.70}{43.74}{56}
\emline{62.70}{43.74}{57}{65}{42}{58}
\put(58,40){\vector(-4,3){0.2}}
\emline{65}{40}{59}{63.35}{38.60}{60}
\emline{63.35}{38.60}{61}{61.74}{38.02}{62}
\emline{61.74}{38.02}{63}{60.15}{38.26}{64}
\emline{60.15}{38.26}{65}{58}{40}{66}
\put(49,14){\vector(4,-3){0.2}}
\emline{42}{14}{67}{43.53}{15.40}{68}
\emline{43.53}{15.40}{69}{45.10}{15.98}{70}
\emline{45.10}{15.98}{71}{46.70}{15.74}{72}
\emline{46.70}{15.74}{73}{49}{14}{74}
\put(42,12){\vector(-4,3){0.2}}
\emline{49}{12}{75}{47.35}{10.60}{76}
\emline{47.35}{10.60}{77}{45.74}{10.02}{78}
\emline{45.74}{10.02}{79}{44.15}{10.26}{80}
\emline{44.15}{10.26}{81}{42}{12}{82}
\put(16,14){\vector(-4,-3){0.2}}
\emline{24}{14}{83}{22.18}{15.40}{84}
\emline{22.18}{15.40}{85}{20.36}{15.98}{86}
\emline{20.36}{15.98}{87}{18.55}{15.74}{88}
\emline{18.55}{15.74}{89}{16}{14}{90}
\put(24,12){\vector(4,3){0.2}}
\emline{16}{12}{91}{17.82}{10.60}{92}
\emline{17.82}{10.60}{93}{19.64}{10.02}{94}
\emline{19.64}{10.02}{95}{21.45}{10.26}{96}
\emline{21.45}{10.26}{97}{24}{12}{98}
\put(68,14){\vector(-4,-3){0.2}}
\emline{76}{14}{99}{74.18}{15.40}{100}
\emline{74.18}{15.40}{101}{72.36}{15.98}{102}
\emline{72.36}{15.98}{103}{70.55}{15.74}{104}
\emline{70.55}{15.74}{105}{68}{14}{106}
\put(76,12){\vector(4,3){0.2}}
\emline{68}{12}{107}{69.82}{10.60}{108}
\emline{69.82}{10.60}{109}{71.64}{10.02}{110}
\emline{71.64}{10.02}{111}{73.45}{10.26}{112}
\emline{73.45}{10.26}{113}{76}{12}{114}
\put(8,41){\vector(1,0){0.2}}
\emline{4}{26}{115}{3.92}{29.64}{116}
\emline{3.92}{29.64}{117}{3.98}{32.78}{118}
\emline{3.98}{32.78}{119}{4.20}{35.41}{120}
\emline{4.20}{35.41}{121}{4.57}{37.55}{122}
\emline{4.57}{37.55}{123}{5.10}{39.18}{124}
\emline{5.10}{39.18}{125}{5.77}{40.32}{126}
\emline{5.77}{40.32}{127}{6.59}{40.95}{128}
\emline{6.59}{40.95}{129}{8}{41}{130}
\put(6,23){\vector(0,-1){0.2}}
\emline{10}{29}{131}{8.09}{28.14}{132}
\emline{8.09}{28.14}{133}{6.79}{26.65}{134}
\emline{6.79}{26.65}{135}{6}{23}{136}
\put(88,26){\vector(0,-1){0.2}}
\emline{84}{41}{137}{85.04}{41.07}{138}
\emline{85.04}{41.07}{139}{85.92}{40.63}{140}
\emline{85.92}{40.63}{141}{86.65}{39.70}{142}
\emline{86.65}{39.70}{143}{87.24}{38.26}{144}
\emline{87.24}{38.26}{145}{87.67}{36.33}{146}
\emline{87.67}{36.33}{147}{87.95}{33.89}{148}
\emline{87.95}{33.89}{149}{88.07}{30.95}{150}
\emline{88.07}{30.95}{151}{88}{26}{152}
\put(82,29){\vector(-2,1){0.2}}
\emline{86}{23}{153}{85.69}{25.47}{154}
\emline{85.69}{25.47}{155}{84.77}{27.32}{156}
\emline{84.77}{27.32}{157}{82}{29}{158}
\put(20,17){\makebox(0,0)[cb]{$u$}}
\put(20,9){\makebox(0,0)[ct]{$\varepsilon$}}
\put(46,9){\makebox(0,0)[ct]{$\delta$}}
\put(46,17){\makebox(0,0)[cb]{$d$}}
\put(72,17){\makebox(0,0)[cb]{$m_r$}}
\put(72,9){\makebox(0,0)[ct]{$\triangle_r$}}
\put(30,45){\makebox(0,0)[cb]{$a$}}
\put(31,37){\makebox(0,0)[ct]{$a^{-1}$}}
\put(62,45){\makebox(0,0)[cb]{$d^*$}}
\put(61,37){\makebox(0,0)[ct]{$\delta^*$}}
\put(8,26){\makebox(0,0)[lc]{$\ev_r$}}
\put(84,27){\makebox(0,0)[rt]{$L$}}
\put(89,31){\makebox(0,0)[lc]{$l$}}
\put(6,42){\makebox(0,0)[rb]{$\coev_l\times\id$}}
\put(46,3){\makebox(0,0)[ct]{\begin{tabular}{c}\textsc{Figure} 18.
Heptagon: from Leibniz to Cartan.\end{tabular}}}\end{picture}\end{center}}

\def\guX{\special{em:linewidth 0.4pt}\unitlength 1mm\linethickness{0.4pt}
\begin{center}\begin{picture}(46,17)      
\put(4,7){\oval(4,4)[b]}
\put(4,11){\oval(4,4)[t]}
\emline{6}{7}{1}{6}{11}{2}
\emline{4}{15}{3}{4}{13}{4}
\emline{4}{3}{5}{4}{5}{6}
\put(2,7){\circle*{1.33}}
\put(2,11){\circle*{1.33}}
\put(42,7){\oval(4,4)[b]}
\put(42,11){\oval(4,4)[t]}
\emline{40}{7}{7}{40}{11}{8}
\emline{42}{3}{9}{42}{5}{10}
\emline{42}{15}{11}{42}{13}{12}
\put(44,11){\circle*{1.33}}
\put(44,7){\circle*{1.33}}
\put(30,10){\vector(2,-1){0.20}}
\emline{16}{10}{13}{18.19}{11.05}{14}
\emline{18.19}{11.05}{15}{20.38}{11.72}{16}
\emline{20.38}{11.72}{17}{22.56}{11.99}{18}
\emline{22.56}{11.99}{19}{24.75}{11.88}{20}
\emline{24.75}{11.88}{21}{26.94}{11.37}{22}
\emline{26.94}{11.37}{23}{30}{10}{24}
\put(0,-2){}
\put(16,8){\vector(-2,1){0.2}}
\emline{30}{8}{25}{27.81}{6.95}{26}
\emline{27.81}{6.95}{27}{25.63}{6.28}{28}
\emline{25.63}{6.28}{29}{23.44}{6.01}{30}
\emline{23.44}{6.01}{31}{21.25}{6.13}{32}
\emline{21.25}{6.13}{33}{19.06}{6.63}{34}
\emline{19.06}{6.63}{35}{16}{8}{36}
\put(23,13){\makebox(0,0)[cb]{$b$}}
\put(23,5){\makebox(0,0)[ct]{$b^{-1}$}}\end{picture}\end{center}}

\def\guY{\special{em:linewidth 0.6pt}\unitlength 1mm\linethickness{0.6pt}
$$\begin{picture}(64,19)         
\put(8,19){\oval(6,12)[b]}
\put(34,19){\oval(6,20)[b]}
\put(60,19){\oval(6,20)[b]}
\emline{8}{13}{1}{8}{5}{2}
\emline{34}{5}{3}{34}{9}{4}
\emline{60}{9}{5}{60}{5}{6}
\put(47,12){\makebox(0,0)[cc]{$+$}}
\put(21,12){\makebox(0,0)[cc]{$=$}}
\put(8,9){\circle*{2}}
\put(31,14){\circle*{2}}
\put(63,14){\circle*{2}}
\put(57,14){\circle{2}}
\put(4,18){\makebox(0,0)[rc]{$A$}}
\put(6,9){\makebox(0,0)[rc]{$D$}}
\put(9,6){\makebox(0,0)[lc]{$A$}}
\put(12,18){\makebox(0,0)[lc]{$A$}}
\put(13,11){\makebox(0,0)[rb]{$m$}}
\put(29,14){\makebox(0,0)[rc]{$D$}}
\put(36,9){\makebox(0,0)[lt]{$m$}}
\put(55,14){\makebox(0,0)[rc]{$\psi$}}
\put(62,9){\makebox(0,0)[lt]{$m$}}
\put(65,14){\makebox(0,0)[lc]{$D$}}
\put(32,3){\makebox(0,0)[ct]{\begin{tabular}{c}\textsc{Figure} 19.
$\psi$-skew algebra derivation $D_A\in\der(m_A)$
\end{tabular}}}\end{picture}$$}

\def\der{\hbox{der}\,}

\def\hom{\hbox{hom}\,}
\newcommand{\N}{\mathbb{N}}

\newcommand{\ie}{\textit{i.e.}$\;$}
\def\zero{\hbox{\bf zero}}
\newcommand{\ra}{\longrightarrow}\def\lra{\longrightarrow}
\newcommand{\id}{\text{id}}      

\newcommand{\cat}{\text{\textbf{cat}}}

\newcommand{\set}{\text{\textbf{set}}}

\newcommand{\obj}{\text{obj}\,}

\def\ev{\hbox{ev}}\def\coev{\hbox{coev}}

\newcommand{\nat}{\text{nat}}

\newcommand{\bt}{\begin{tabular}{c}}
\newcommand{\et}{\end{tabular}}
\newcommand{\ba}{\begin{array}{c}}
\newcommand{\ea}{\end{array}}
\def\bt{\begin{tabular}{c}}\def\et{\end{tabular}}

\begin{document}\title{Categorical Analysis}\author{Zbigniew Oziewicz}
\address{Universidad Nacional Aut\'onoma de M\'exico, Facultad de
Estudios Superiores Cuautitl\'an, Apartado Postal \#25,
C.P. 54700 Cuautitl\'an Izcalli, Estado de M\'exico,
and Uniwersytet Wroc{\l}awski, Poland}\email{oziewicz@servidor.unam.mx}
\author{Guillermo Arnulfo V{\'a}zquez Couti{\~n}o}
\address{Universidad Aut\'onoma Metropolitana, Unidad Iztapalapa,
Avenida Michoac\'an y la Pur\'{\i}sima, Colonia Vicentina, Apartado
Postal 55-534, C.P. 09340 M\'exico D.F.}\email{gavc@xanum.uam.mx}
\date{December 7, 1999}
\thanks{Supported by el Consejo Nacional de Ciencia y Tecnolog\'{\i}a
(CONACyT) de M\'exico, grant \# 27670 E, and by UNAM, DGAPA, Programa de
Apoyo a Proyectos de Investigaci\'on e Innovaci\'on Tecnol\'ogica,
Proyecto IN-109599 (1999-2002).\\Zbigniew Oziewicz is a member of Sistema
Nacional de Investigadores, M\'exico, No. de expediente 15337.}

\hyphenation{ca-te-go-ry Wo-ro-no-wicz non-co-m-mu-ta-tive}

\begin{abstract} We propose the categorification of algebraic analysis
in terms of a specific 2-category, called here the Leibniz 2-category
given by generators and relations which include the Leibniz-like
relation (strict 3-cell) among extended 2-cells. The Leibniz 2-category
offers the `most general' notion of a `(co)derivation', as a strict
3-cell, for a general (al- co-)gebra, not necessarily (co-)associative,
not necessarily (co-)unital, nor necessarily (co-)\-co\-m\-mu\-ta\-tive.

We outline a program in which every 2-cell related to a partial
(co-)derivation (called a Leibniz strict 3-cell) is translated into an
appropriate 2-cell related to a Cartan's-like (co-)derivation (called a
Cartan strict 3-cell), and vice versa. We found also that a
non-Leibniz component (2-cell related to a Leibniz 3-cell), responsible
for the stochastic calculus, must be a kind of ternary operation.
\end{abstract}\maketitle

\swapnumbers
     \theoremstyle{definition}
\newtheorem{example0}{Pedagogical Example}[section]
\newtheorem{rem1}[example0]{More History}
\newtheorem{def1}{Definition}[section]
\newtheorem{example1}{Example}[section]
\newtheorem{example2}[example1]{Example}
\newtheorem{def2}[example1]{Definition}
\newtheorem{def3}{Definition}[section]
\newtheorem{note1}[def3]{Note}
\newtheorem{bigebra}[def3]{Bigebra}
\newtheorem{motivation}[def3]{Motivation}
\newtheorem{geometry}{Geometry}[section]
\newtheorem{problem}{Problem}[section]
\newtheorem{rem2}[problem]{Remark}
\newtheorem{def4}[problem]{Definition}    
   \theoremstyle{plain}
\newtheorem{Tarski}[example1]{Theorem} 
     \theoremstyle{definition}

\tableofcontents

\section{Introduction} A general algebra $A,$ not necessarily unital,
associative nor commutative, with one partial derivation
$D\in\der(A,A)\subset A^A$ is said to be a differential algebra [Ritt
1950] or a Leibniz algebra.

An algebra $A$ together with a set (eventually a Lie algebra) $L$
acting as partial derivations on an algebra $A,$
$L\rightarrow\der(A,A),$ do possess even a richer historical
terminology showing the importance of this concept in algebraic
analysis, and in commutative and non-commutative differential geometry.
Here is the partial list of a multitude of names, however, the
corresponding references are mostly omitted because this historical
remark is irrelevant for what follows. A pair of algebras $(A,L)$ is
said to be a pseudo-alg\`ebre de Lie [Herz 1953], a differential Lie
algebra [Palais 1961], an $A$-Lie algebra [Rinehart 1963], a Lie module
[Nelson 1967], a Lie-Cartan pair [Kastler \& Stora 1985], a Leibniz pair
[Flato, Gerstenhaber and Voronov 1995], a Cartan pair [Borowiec 1996], a
Lie-Rinehart algebra [Huebschmann 1999].

What is known as the Leibniz `rule', or as the Leibniz condition
defining a partial or a Cartan's derivation, we prefer to call rather as
the Leibniz relation, Leibniz's axiom, in analogy to the presentation in
terms of `generators and relations' in universal and linear algebras.
Algebraic analysis deals exactly with Leibniz-like relations,
\ie with various deviations from the strict Leibniz relation in terms of
the quasi-Leibniz and non-Leibniz components, or such deviation as there
is in a stochastic calculus [Arnold 1974, \S 5.3-5.4, pp. 89-91; Sobczyk
1991]. It is not our intention,
however, to go into stochastic calculus in the present paper. In the
standard set theory, if $A$ is a binary algebra and $x,y\in A,$ then the
Leibniz relation for a partial quasi derivation $D\in q\der(A,A)$
is presented in terms of irrelevant elements in the following form
\begin{equation}\forall\,x,y\in A,\qquad D(xy)=(Dx)y+xDy+(Dx)(Dy)+\ldots.
\label{Leib}\end{equation}
Because this holds for all elements of $A$ and for any quasi derivation
in $q\der(A,A),$ therefore the three symbols, $x,y$ and $D$
in (\ref{Leib}), are in fact totally irrelevant for this relation. If we
drop them what is left? The genuine Leibniz-like relation (among what
species, quantities?) is hidden in expression (\ref{Leib}). Even more,
an entire algebra $A$ must be considered as a variable object (and not
only the elements of $A$) because the Leibniz-like relation holds
probably for many algebras. We wish to express the Leibniz relation not
only as element free, but also as independent of the choice of an
algebra, as object-free, \ie as a model independent abstract relation,
exhibiting exactly what this relation is, and among what species.
We see analysis as a part of universal algebra, rather than a part of
just linear algebra. If we wish to consider an algebra (a ring) as an
irrelevant variable object then we must move from functions and
operations to categories and functors. Essentially, we do propose the
categorification of analysis, a categorical analog of algebraic
analysis, \ie the categorical analysis. A guiding motto for the passing
XX century by Gian Carlo Rota [1998, p. 3] `analysis play second fiddle
to algebra', must be supplemented for the next XXI century with `algebra
play second fiddle to categories'. Categories are an aid to understanding.

\begin{example0} We wish to illustrate the general
idea on an example of the bilinear term $xy\in A\in\obj\cat$ in
(\ref{Leib}). First, $xy$ is the value of the element-free
multiplication $m_A:A\otimes A\rightarrow A,$ \ie $m_A(x\otimes y)=xy.$
Next, we do not need to be restricted to some specific bifunctor
$\otimes:\cat\times\cat\rightarrow\cat,$ therefore we will use a
symbol $\guy$ denoting the name of an arbitrary bifunctor (with a place
on the top for {\it two} inputs and a place at the root for {\it one}
output). In an expression for a value of a multiplication $m_A,$ a
bifunctor $\guy$ is restricted to (or composed with) a diagonal
cofunctor $\guc:A\mapsto(A,A),$ and this composition (we prefer the
name {\it grafting}) gives a composed endofunctor denoted pictorialy by
$\guo\equiv\guy\circ\guc$. Therefore $m_A$ is a map from an object $\guo
A$ to an object $(\id)A$ and finally an object-free multiplication $m$
(which we wish that substitute the $xy$ term in (\ref{Leib})) is a
(natural) transformation from an endofunctor $\guo$ to an identity
functor $|\equiv\id\equiv\id_{\cat}.$\end{example0}

\begin{rem1}[More History] Hausdorff proposed for noncommutative
algebra to keep the Leibniz axiom still in the same form (\ref{Leib}) as
for a commutative algebra. Gian-Carlo Rota, Sagan and Stein found in
1980 that the Leibniz axiom for non-commutative algebra in the
Hausdorff form (\ref{Leib}) break down the chain rule for partial
derivation (derivative), and they proposed an altogether different
notion of derivative, which is {\it not} a derivation, called the {\it
cyclic} derivative with the cyclic axiom instead of the Leibniz axiom
(\ref{Leib}). The cyclic axiom implies that the chain rule still holds
for non-commutative algebra.

Altogether different partial derivations for a not necessarily
commutative algebra were proposed by Woronowicz [1989], Majid [993, 1995
\S 10.4], Oziewicz, Paal \& R\'o\.za\'nski [1995], Borowiec [1996,
1997], see {\sc Motivation} 5.4 below and the last Section 12.\end{rem1}

Our purpose is to describe some examples of the 2-categories given by
generating 1-cells and by generating 2-cells. A Leibniz 2-category
includes the Leibniz-like binary relation (axiom) on expanded 2-cells,
and such axiom is nothing more than a strict 3-cell. Still weaker is a
Leibniz 3-category which includes the not strict 3-cells of the
modifications of the Leibniz axiom. Another purpose of this paper is to
suggest some directions for the study of 2-and 3-categories with
Leibniz-like axioms.


\section{2-Category} The classical works on 2-categories are
[B\'enabou 1967, Gray 1974; Kelly \& Street 1974]. For a recent account
we refer to [Baez \& Dolan 1998, Batanin 1998].

A set of natural numbers with zero (non negative integers) is denoted
by $\N\ni 0.$ For $n\in\N,$ let $G_n$ be a collection of `$n$-cells',
such that $G_{n+1}$ is a disjoint sum of sets indexed by $G_n\times
G_n.$ If $x,y\in G_n\times G_n$ then the corresponding subset in
$G_{n+1}$ is denoted by $G_n(x,y)\subset G_{n+1}.$ Let $B{\guv}A$ be a
convenient notation for a `bundle' $(B\times B)\leftarrow A.$ An n-graph
is $G_0{\guv}G_1{\guv}G_2{\guv} \ldots{\guv}G_n,$ etc, and similarly for
an $\infty$-graph.

\begin{def1} An n-category is an n-graph with (strictly or
weakly) associative and unital compositions for all
$i\in\{0,1,\ldots,n-1\},$
\begin{equation}\hbox{for all}\;x,y,z\in G_i,\qquad G_i(x,y)\times
G_i(y,z)\lra G_i(x,z).\label{comp}\end{equation}\end{def1}

This means that for each fixed $x,y\in G_i\times G_i,$ a collection
$G_i(x,y)\subset G_{i+1}$ with
$$\forall\alpha,\beta\in G_i(x,y)\subset G_{i+1},\qquad
[G_i(x,y)](\alpha,\beta)\;\subset\;G_{i+2}$$ is a 1-category
$\{\hbox{objects}=G_i(x,y),\hbox{morphisms}=\ldots\}.$ In particular
$$\cat\equiv\{\hbox{objects}=G_0{\guv}G_1=\hbox{morphisms}\}$$ is a
1-category and a 2-category $G_{0,1,2}$ is a collection of 1-categories
$\{G_0(x,y)\}$ indexed by $G_0\times G_0.$

For the case $G_i\equiv\N,$ the only one considered in this paper, we
found convenient for an (i+1)-cell to adopt the name {\it plant} from
botany, because plants will be {\it grafted}, and not only composed
(\ref{comp}), (see Section 4). An (i+2)-cell is called a `formal
instance' in [Kelly \& Laplaza 1980, \S 4, p. 199], or an arrow, or the
{\it name} of a natural transformation, because models (with respect to
n-functors) of the (i+2)-cells are natural transformations of
(i+1)-cells. In this paper, for an (i+2)-cell we adopt the name
{\it footpath} (among plants).

A 2-category of plants and his footpaths, with grafting of plants,
but with the usual (strict) associative compositions of footpaths
(\ref{comp}), can be called {\it an algebraic 2-category with grafting}
$\{\N,\hbox{plants},\hbox{footpaths}\}.$ However, we adopt the shorter
name `grafted club', where the name `club' was introduced by Kelly
[1972, pp. 113--123], although Kelly's formal definition of a club is
not exactly the same as our's.

\section{Plants} If for some $i\in\N,$ a collection of i-cells
$G_i$ has an extra structure of `an algebraic sketch', \ie is a
(strict or weak) monoidal category whose monoid of object is $\N,$
with $+$ for the binary monoidal operation and with terminal $0\in\N,$
then an $(i+1)$-cell is said to be {\it a plant,} and the collection of
all plants $P\equiv G_{i+1}$ is a disjoint collection of 1-categories
indexed by $\N\times\N.$ Every such 1-category is denoted by
$(m\mapsto n),$ $$\obj(m\mapsto n)\subset P.$$
In this case, a monoid $\N$
with a collection of (i+1)-morphisms $G_{i+1}\equiv\{n\mapsto m\},$
where $n\mapsto m$ is a 1-category, is a monoidal small 2-category.
Every plant is an object of a unique 1-category $m\mapsto n$ for some
$m,n\in\N\times\N$ and denoted by (Figure 1 suppose strict associativity
of $\N\times\N\stackrel{+}{\rightarrow}\N$)

\guA   

All plants (as graphs) are implicitely directed from the top (input:
leaves from $\N$) to the bottom (output: roots from $\N$). Plants do
not have outer nodes, every node is inner. One can think that a plant
$\in P$ is the {\it name} of a functor, because models of plants are
many variable and many-valued multi-functors often of mixed variances.
Our plant is almost the same as a `shape' introduced by Kelly and
MacLane [1971, p. 102].

In this paper we can rather restrict ourself to the case of one sort (or
one color), \ie in the realization to just one (but any) category $\cat$
with structure (one exception will be discussed separately late on).
In such cases of one sort, the edges need not to be colored, and in a
model category $\cat,$ all edges of plants are the same identity functor
$\id\equiv{\id}_{\cat}.$ A plant $\id\in(1\mapsto 1)$ given by a
vertical dash on Figures 2 \& 4 (or any other edge without vertices)
in a model category $\cat$ is mapped to $\id\equiv{\id}_{\cat}.$

A plant on Figure 1, given by a directed tree with vertices, with $m$
free leaves and $n$ free roots, \ie an element of a category $n\mapsto
m,$ in a model is mapped to a multi-variable and multi-valued functor
(or morphism): input $\rightarrow$ output, $\cat^{\times
n}\rightarrow\cat^{\times m}.$ In particular, any endofunctor on $\cat,$
different from $\id,$ like the adjoint in a compact closed category does
represent a plant from $1\mapsto 1$ with exactly one input and one
output and with exactly one vertex (labelled if needed) in the case of
not being composed.

A category with one object and one (identity) morphism is said to be
the zero category $\zero\equiv\cat^{\times 0}$ (some authors do prefer
the name unit category, as in [Kelly 1972, p. 74]). For example, if
$\Bbbk$ is a field, then in the category $\Bbbk$-{\bf vec} of
the vector $\Bbbk$-spaces, we have ($\Bbbk$-{\bf vec})$^{\times 0}\equiv\Bbbk,$
whereas in the category of sets {\bf set}$^{\times 0}\equiv\emptyset.$
A category $1\mapsto 0$ possesses the unique plant which we call
{\it killer.} The killer in a model is represented by the functor from
$\cat$ to the terminal $\zero.$

\guB  

A plant from the category $0\mapsto 1$ we call a {\it creator}. In a
model, a creator is represented by a functor from initial $\zero$ to
$\cat$ and is identified with some object of $\cat.$ These plants are
given on Figure 2 by edges with one vertex, they are objects (elements)
of categories $1\mapsto 0$ and $0\mapsto 1$ respectively, what for the
reader convenience is indicated on the bottom of the plants. The same
convention of indicating a category, to which a given plant belongs, at
the bottom of a plant, is used on Figures 3 \& 4.

A plant on the right in the second row on Figure 2 in realization is a
bifunctor $\cat^{\times 2}\rightarrow\cat,$ and a plant on the left in
realization is a `cobifunctor' $\cat\rightarrow\cat^{\times 2}.$
We wish to have a functorial calculus with explicit functors whose
co-domain is a product $\cat\times\cat$, like the diagonal functor
$\cat\rightarrow\cat\times\cat$ given by $\cat\ni S\mapsto(S,S).$
A functor $\cat^{\times m}\rightarrow\cat^{\times n}$ is the same as
an n-tuple of functors $\cat^{\times m}\rightarrow\cat.$ A
functor $\cat\rightarrow\cat^{\times n}$ is an ordered set (of
cardinality $n$) of $n$ endofunctors, therefore there is a bijection of
functor categories,
$$(\cat_A\times\cat_A)^{\cat_B}\quad\simeq\quad
(\cat_A)^{(\cat_B+\cat_B)}.$$
Here $+$ is the disjoint sum (coproduct). This means that we can get by
with a calculus that has no explicit place for a disjoint sum of
categories - and that is precisely what we wish.

\section{Grafting} In the case, when $G_i=\N$ is a monoid, we wish
to generalize the composition (\ref{comp}) to a multivalued operation
for which we adopt the name {\it grafting} from botany. Consequently,
the objects of a 1-category $n\mapsto m$ ((i+1)-cells, (i+1)-morphisms),
we call {\it plants} (processes, etc), because plants are going to be
allowed to be `composed' in a more sophisticated way than compositions
of morphisms. Let $P$ be a collection of plants, \ie a collection of
1-categories $\{m\mapsto n\}$ indexed by $\N\times\N.$ A multivalued
operation of grafting generalize and unify composition, substitution and
concatenation. The grafting of plants is a specific {\it symmetric} map,
\begin{equation}\hbox{grafting in an orchard:}\quad P\times P\quad\lra
\quad 2^P.\label{grafting}\end{equation}
The above map generates an expanded grafting, \ie an associative,
unital and commutative binary operation $2^P\times 2^P\lra 2^P.$
One can graft {\it any} two plants and they do not need to be
necessarily of the same arities. Grafting a plant from a category
$n\mapsto m$ with a plant from $p\mapsto q$ (with no restrictions like
$m=p$ and so on at all!) gives no more than $n+p+m+q$ different derived
plants. However, these derived plants are scattered among different
categories $\{n+p\geq r\mapsto s\leq m+q\}.$ The grafting is a
collection of $n$-graftings, say, components of grafting. However, in
general, any specific $n$-grafting can be still multivalued. In
particular a $0$-grafting is the same as a 2-valued concatenation.

Instead of a formal definition of grafting we will explain what we have
in mind on two examples.

\begin{example1} Grafting a plant $\guy\in (2\mapsto 1)$ with
a plant $\guc\in(1\mapsto 2)$ (in any order, because grafting is the
commutative `operation') gives six derived plants as shown on Figure 3.
All these derived plants are two-letter words without repetition of the
same letter twicely.

\guC  

To get the above six derived plants one must `rotate and graft'
one plant around another (without rotating the plants itself) up to
isotopy (of the kind of Reidemeister moves). Therefore, concatenation
of plants giving two composed plants from $3\mapsto 3$ on Figure 3 (the
only not connected graphs or $0$-grafted), can be considered as a
particular case of a grafting, namely this is a $0$-grafting. In Figure
3 we have three composed plants of type $2\mapsto 2$ derived by
$1$-grafting and one composed plant from $1\mapsto 1$ derived by
$2$-grafting. The derived plants on Figure 3 are ordered according to
rotation.\end{example1}

\begin{example2} The grafting of a plant from $1\mapsto 1$ with
itself gives only two new plants, one from $1\mapsto 1$ (composition =
1-grafting) and another from $2\mapsto 2$ (concatenation =
0-grafting).\end{example2}

The (free) generators of all plants are said to be fundamental (or
basic) plants (alphabet of plants), and similarly the generators of all
footpaths are called fundamental or basic footpaths (alphabet of
footpaths). An alphabet of plants generates a free club of derived
plants, derived 1-cells [Kelly 1972, \S 3, p. 116, 1974] (or an operad,
garden, forest, orchard), made from fundamental plants by iterated
grafting, generalizing and unifying substitutions, compositions and
concatenations.

An alphabet of plants does freely generate by iterated graftings
an infinite garden-forest of derived plants-words ($\N$-graded by
the number of letters, $\N$-graded by the number of graftings and
$\N\times\N$-graded by input \& output). The set of all plants derived
in this way, with the multi-valued operation of grafting, is said to be
a free club (or operad) presented by this alphabet. (We feel that this
formalism needs also co-grafting, a split plant $2^P\rightarrow
2^P\times 2^P$; however, this is another story.)

\begin{def2} An $n$-category $G$ (or an $\infty$-category
$G$) is said to be a {\it plant-like} if the following two conditions
holds:
\begin{description}
\item[(i)] $\;$For some $i\in\N,$ $G_i\equiv\N,$ \ie
$$G\equiv\{\ldots,{\guv}G_{i-1}{\guv}\N{\guv}(P\equiv G_{i+1}){\guv}
G_{i+2}{\guv}\ldots\}.$$
\item[(ii)] The usual partial compositions (\ref{comp}) for plants
$P\equiv G_{i+1}$ are extended to multivalued grafting $P\times P\lra
2^P$ (\ref{grafting}).\end{description}
If $G$ is a plant-like $n$-category with $G_i\equiv\N,$ then a
$(i+1)$-cell in $P\equiv G_{i+1}$ is said to be a {\it plant.}\end{def2}

\section{Why?} Why [should one] use graphs (tangles, hieroglyphics)
jointly with the techniques of naming objects by letters from the Greek
or Latin alphabets?

Graphs convey more information than these letters. The techniques of
naming objects and arrows (elements, functors, morphisms, functions,
etc) do allow a large amount of the most essential but routine detail
to be hidden, like changing street's names after a political revolution.

Graphs lead to drawing pictures as in this paper, which display
the relationships between various operations, functors etc. This must be
contrasted with the usual lists of Greek or English letters and
equations which frequently convey nothing at all to the reader. Compare
for example graphs from Figure 11 below, with a Table of the letters
after Motivation 6.4. If there is nothing else, we do like graphs and it
is easier to spot errors in a graph than in a list. We agree with
Dieudonn\'e: `half the success (in mathematics) depends on a proper
choice of notation'. Unfortunate notation might kill fortunate ideas.

\section{The Leibniz alphabets: plants and footpaths} In what follows we
will need an alphabet consisting of four plants: $\id\in(1\mapsto 1)$,
a killer, a binary tree and a two-rooted tree as shown on Figure 4.
It is convenient to call this alphabet the Leibniz alphabet - this
choice is motivated at the end of this Section.

\guD  

The Leibniz alphabet of plants given by Figure 4, generates by
iterated graftings a free algebra (club) on four plants, the garden of
derived plants (terms), and a small fragment of this garden, relevant
for what follows, namely the fragment of type $2\mapsto 1$ is shown
next

\guE  

Among the seven plant-words on Figure 5, five of them are grafted from
three letters.

\begin{def3} If $G$ is plant-like $n$-category with
$G_i\equiv\N,$ then an $(i+2)$-cell $\in G_{i+2}$ is said to be a
{\it footpath.} Related names: formal instance or arrow [Kelly \&
Laplaza 1980, p. 199]. A footpath is a map which does not change the
arity of a plant,$$(m\mapsto n)\ni\hbox{plant}_1\;\mapsto\;
\hbox{plant}_2\in(m\mapsto n).$$ A not invertible footpath is said to
be `lax', a twosided invertible footpath also is said to be a `pseudo
relation', and a strict (equivalence) relation is the same as an
identity $\approx.$\end{def3}

In what follows we use the convention that not named paths are
necessarily twosided invertible. A model of a footpath is a natural
transformation of functors, more about models is in Section 11.

We need an alphabet of footpaths. We wish first to discuss the
following footpaths as a fundamental one (a motivation for this
choice will be explained later on),

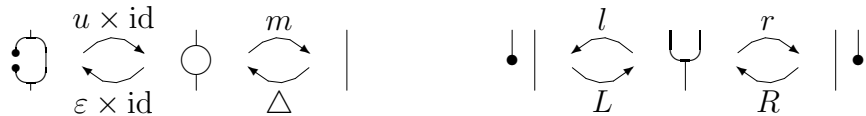
\guF  

The above footpaths we call appropriately: a binary algebra $m$ (a
multiplication), a binary cogebra $\triangle$ (a comultiplication), the
left and the right actions $l$ \& $r,$ the left and the right
co-actions $L$ \& $R.$ A path $u$ is said to be a left unit for an
algebra $m$ (or an algebra $m$ is said to be a left $u$-unital) if a
composed 2-cell $m\circ(u\times\id)$ is twosided invertible.
We do not yet wish to impose this pair of two strict 3-cells-axioms,
because exists the competing possibility of a weaker left unit in the
case that this composed 2-cell is only one-sided invertible. We are
going to motivate these names next. Note that, disregarding
universality, $l$ \& $r$ looks like a product, \ie as a pair of
`projections'. Similarly $L$ \& $R$ looks like a coproduct \ie a sum.

We will explain an interpretation of these footpaths on an example of
the particular functorial realization of the plants $\guy$ and $\guc.$
This example will give a feeling, what these footpaths could mean in
the usual life of.

A plant $\guy$ can be realized as some bifunctor $\cat^{\times 2}\lra
\cat.$ For example in a category of sets $\set,$ for $A,B\in\obj\set,$
we could have for example
\begin{description}
\item[-] a cartesian product, $\guy(A,B)=A\times B.$\\In a category of
bimodules this is the bifunctor of the tensor product $\otimes.$
\item[-] an exponentiation, $\guy(A,B)=B^A,$\\or more generally
$\guy=\hom$ bifunctor, $\guy(A,B)=\hom(A,B).$
\item[-] any derived bifunctor as $\guy(A,B)=2^A\times B,$ etc.
\end{description}
Let for simplicity $\guc$ be the diagonal, $\guc:\obj\cat\ni
S\mapsto\guc(S)=(S,S).$ A 2-grafting $\guc$ with $\guy$ gives the unique
endoplant $\guo:\cat\lra\cat.$ For a choice in a $\set,$ $\guy=\times,$
an endoplant $\guo$ is given by the value $\guo(S)=S\times S$
(the obvious values on morphisms are omitted). The value of the plant
`$|\equiv\id$' on $S\in\obj\cat,$ (or the evaluation of $S$ on the
plant `$|$') must be $S.$ Therefore, in this case, the value of $S$ on
a footpath $m$ (or vice versa) is a map $m_S:S\times S\rightarrow S,$
\ie $m_S$ is a binary multiplication. A domain of a footpath $m$ is
said to be a type of algebra $(S,m_S),$ and in the above example, this
is a binary algebra.

\begin{note1} A footpath can act on a grafted plant {\it independenly}
of his action on an alphabet of plants. In the example
above, a footpath $m$ do not change neither $\guy$ nor $\guc$ whereas is
transforming their composition $\guo.$\end{note1}

Besides of the above six fundamental arrows $\{m,\triangle,l,L,r,R\},$
we will need in what follows, also two more footpaths, not necessarily
invertible, as shown on Figure 6,

\guG   

An expansion on Figure 6 (and further on Figures 8-11, 13, 15 etc.) is
a particular case of a pasting introduced by B\'enabou in 1967.

The footpaths on Figure 6 might be consequences of the following
stronger conditions: the strict mitosis for the diagonal functor
(Figure 8) and the strict associativity for a bifunctor \guy.
The lax associativity means commutativity of the MacLane pentagon on
Figure 9, where a footpath `a' is stronger than that on Figure 6.
However, in what follows we will need explicitly only the footpaths from
Figure 6.

\guH   

\guI   

All arrows are expandible to the appropriate derived plants.
An alphabet of footpaths provide generators (of a partial monoid)
of expanded footpaths (`multi-arrows', expanded arrows). This includes
in particular a pasting introduced by B\'enabou [1967]. Footpaths need
not to be always composable.

Grafting on plants and the pasting and the expansion on footpaths is a
functor from a 2-category of alphabets to a free 2-category,
$$\{\N,\hbox{alphabets of plants \& footpaths}\}\quad\lra\quad
\{\N,\hbox{derived plants \& pasted footpaths}\},$$
\ie 0-cells are the same, but 1-cells are to be derived and 2-cells are
to be expanded.

\begin{bigebra}[Bigebra] A famous expansion of an associativity path
$a$ provide the MacLane and Stasheff pentagon [MacLane 1963].

\guJ  

A particular example of an expansion of an algebra and cogebra paths
$m$ \& $\triangle$ is shown on Figure 10. This expansion with an
appropriate axiom for 2-cells is said to be a $\sigma$-braided bigebra.

\guK  
\end{bigebra}

The fundamental footpaths $\{m,l,r\}$ can be expanded as follows (an
expanded footpath from Figure 6 is not shown for simplicity on Figure
11 where the `heptagon' must possess 11 edges),

\guL   

\unitlength=1mm

\begin{motivation}[Motivation] The above choices of the Leibniz alphabets
of plants and of footpaths have the following motivation.

Borowiec [1996, 1997] proposed an altogether different notion of a
partial derivation for a not necessarily commutative algebra. Borowiec
postulated the classical Leibniz strict 3-cell-axiom for the Cartan
bimodule-valued derivation $d\in\der(m,m_l\,\&\,m_r)$ and for a
bi-co-module-sourced co-derivation
$\delta\in\der(\triangle_l\,\&\,\triangle_r,\triangle)$
[Borowiec \& V\'azquez Couti\~no], in the framework of the calculus of
the differential forms for a not neccesarily (co- bi-)\-commutative
(al- co- bi-)gebra $m$ \& $\triangle$ [Woronowicz 1989]. Then Borowiec
derived the `correct', but still an altogether different notion of a
partial derivation (which in our terminology is an example of a
3-cell), such that in the case of a commutative algebra the Borowiec
partial derivation coincides with the classical Leibniz partial
derivation.

Please note that the third term in the expression (\ref{Leib}) needs a
permutation of letters $D$ with $x,$ (or sometimes a cyclic permutatiom
as in [Przeworska-Rolewicz 1995 [36], p. 779, formula (3.9), where
$c$ must read $x$]. This is reflected in the functorial calculus in
the necessity of introducing a braid plant - functor
$\gub:\cat^{\times 2}\rightarrow\cat^{\times 2}.$ Contrary to this, the
Borowiec partial derivation do not need permutation, at the cost that in
the first two terms in (\ref{Leib}) $D(xy),Dx\in A$ as usual; however,
the third term, contrary to the habit, must be $(Dx)y$ with $Dx\in L.$
Then Borowiec is proving that the classical Leibniz relation
(\ref{Leib}) appears to be the consequence of the commutativity in an
algebra. The quasi-Leibniz term in (\ref{Leib}) needs besides the braid
$\gub,$ and a $D$-dependent family of multiplications, also duplication
of the letter $D.$ In order to understand the Borowiec results we do not
need braid plant $\gub$ and in this paper we restrict the attention to
the Leibniz-like axioms which do not involve $\gub.$\end{motivation}

In order to categorificate the Leibniz-like 3-terms relation for a
partial derivation as invented by Borowiec [1996, 1997] and
investigated in [Borowiec \& V\'azquez Couti\~no] for a not necessarily
commutative algebra, it appears that we need exactly an alphabet of
plants given by Figure 4, and the fundamental paths as presented above
in this Section.

Altogether we have the following seven paths from plant $\gua\;$ to
plant $\gur\;$, and one must read them from left to right,

\begin{center}\begin{tabular}{cc}\hline
Leibniz's paths&\bt $l\circ(\id\times m)\circ a$\\
$m\circ(l\times\id)$\\$c\circ l\circ(l\times\id)$\et\\\hline
Borowiec's paths&\bt$l\circ(\id\times c)\circ(r\times\id)$\\
$c\circ l\circ(r\times\id)$\et\\\hline
Stochastic paths&\bt$m\circ l\circ a$\\
$c\circ l\circ l\circ a$\et\\\hline\end{tabular}\end{center}

On Figure 11 there are two paths, through the dashed vector, involving
the footpath $l$ represented as a ternary operation `of a type
$3\mapsto 2$'. We call them {\it the stochastic paths} (or non-Leibniz
paths), because they are absent in the standard Leibniz axiom and
because they are responsible for the stochastic differential calculus
[Arnold 1974, \S 5.3-5.4, pp. 89-91; Sobczyk 1991].

In what follows, for simplicity, we identify on a plant $\guo,$ the
footpath $m$ with $c\circ r,$ see Figure 7. The footpaths from Figure
11 in a model category are represented as operations whose graphs are
shown on Figure 12.

\unitlength=1mm
\guM   

\section{The Cartan club} An extended algebra $m\mapsto
m_l\,\&\,m_r,$ is said to be a $m$-bimodule. An extended cogebra
$\triangle\mapsto\triangle_l\,\&\,\triangle_r$ is said to be a
$\triangle$-bicomodule. An extension of a 2-cell is given by a set of
3-cells (axioms if they are strict) [Eilenberg 1948, Gugenheim 1962,
Kelly 1972, p. 94, Cuntz \& Quillen 1995]. A Cartan's-like club (the
Cartan differential calculus of the differential forms) consists,
roughly speaking, of a `derivation' $d\in\der(m,m_l\,\&\,m_r)$ with
target in a $m$-bimodule $m_r$ \& $m_l$ , or/and of a `coderivation'
$\delta\in\der(\triangle_l\,\&\,\triangle_r,\triangle)$ with source in
a $\triangle$-bicomodule $\triangle_r$ \& $\triangle_l$.

Here we use the symbol `der' to denote a Leibniz-like (quasi) axiom for
{\it both} directions, either for a `derivation' as well as for a
`coderivation', although Borowiec \& V\'azquez Couti\~no proposed a
longer symbol `coder'.

The following statements are equivalent:
\begin{eqnarray*}
\hbox{(i)}&d\in\der(m,m_l\,\&\,m_r)\\
\hbox{(ii)}&\hbox{the set of 2-cells}\;\{d,m,m_l,m_r\}\;\hbox{is
related by a 3-cell `der'.}\end{eqnarray*}
A 3-cell `der' is just a particular Leibniz-like axiom. A
categorification of analysis is an interpretation of `der' as a 3-cell
(or as a set of the 3-cells), \ie as a modification(s) in the
terminology of B\'enabou [1967].

Limited space does not allow to discuss 3-cells
in this paper. Therefore our purpose is limited to the
discussion of 2-cells only.

The Cartan club does need, besides the Leibniz alphabet given on
Figure 4, just one more plant which must be a creator. Altogether, we
need five fundamental plants. The most important fragment of the derived
garden of type $1\mapsto 1,$ what we call the Cartan garden, is shown on
Figure 13, together with the all involved now fundamental footpaths
$\{d,m,m_l,m_r,\ldots\}.$

A fundamental footpath denoted by $d,$ subject to the appropriate
axiom, is going to be the Cartan bimodule-valued derivation with
respect to the triple $\{m,m_r,m_l\}.$ The Cartan footpath $d$ looks
like a `gluing slot', however the right interpretation must be that $d$
is `transforming' a constant functor into the identity functor.

In the Cartan garden on Figure 13 we have, among other, three classical
Cartan's paths, and also a path which we call `the stochastic
differential' - this is non trivially expanded from the Cartan
fundamental footpath $d.$

\guN  

The two plants are named sink and source, and these names correspond to
a `derivation' arrow $d.$ For a coderivation $\delta$ these names must
be interchanged: sink $\leftrightarrow$ source.

\subsection{Adjoint plant} The Cartan-Leibniz and
Leibniz-Cartan dictionaries, which will be given in the next Sections,
need two extra fundamental plants, which, strictly speaking belongs
neither to the Leibniz garden, nor to the Cartan garden, however must be
in a dictionary. The first extra plant is $\guj\in(1\mapsto 1),$ to be
refered as an adjoint with respect to a unital binary plant and some
suitable footpaths.

A unit \& counit for the binary plants \guy \& \guc is the same as an
invertible (co-)modul paths, as for example the appropriate 2-cells on
Figure 13 in case they are invertible. We keep the convention that not
named paths are necessarily twosided invertible, as for example the
invertible paths on `cyclic pentagons' with the expanded evaluation \&
coevaluation on Figure 15.

The following definition was essentially invented by Kelly \&
La\-pla\-za [1980].

Let a creator $\gcr\in(0\mapsto 1)$ be a unit for a binary plant \guy.
A plant $\guj\in(1\mapsto 1)$ is said to be {\it left adjoint} with
respect to a unital \guy \& \gcr, given jointly with \guc, if there
exist a pair of fundamental footpaths called left evaluation $\ev_l,$
and a left coevaluation $\coev_l$, such as shown on Figure 14.

\guO  

The two fragments of the garden $1\mapsto 1$ on Figures 14 \& 15 with
evaluation and coevaluation, and with not named paths being invertible
unit \& counit is known as the closed structure [Kelly \& Laplaza 1980,
p. 193]. For evaluation \& coevaluation paths is necessary that a
binary plant \guy\, be unital.

In [Kelly \& Laplaza 1980] the evaluation is called a counit and denoted
by $e$, and the coevaluation is called a unit and denoted by $d$ (what
can not be accepted here because we wish to keep the symbol $d$ for a
footpath which under an additional axiom play the role of the
traditional bimodule-valued Cartan derivation). However, a twosided
`counit' for a plant \guc is a pair of identities given by mitosis on
Figure 8.

\guP  

\begin{geometry}[Geometry] In the most interesting applications in
functorial realization appears to be that plants \guj and
$|\in(1\mapsto 1)$ do have different variance. In this case, a path
between them must be modeled by, what is called, a Barr dinatural
transformation [Par\'e \& Rom\'an 1998] (see Section 10 for more
details). We use for this footpath, on Figures 15 \& 16, the name
`geometry', because in a particular important example of different
variance, the geometry footpath is responsable for a Riemann-like
structure (a scalar-like product) in the riemannian or euclidean
geometries. A geometry footpath would be necessary if we would like to
have the concepts of gradient and rotation in the differential
calculus. However, the riemannian differential geometry is outside the
scope of the present paper and a geometry footpath will not be
considered in what follows.\end{geometry}

An unipotent plant \guj (Figure 2) is also said to be reflexive. When
drawing the right evaluation pentagon (Figure 15) we assumed
implicitely, just for simplicity, that an adjoint plant \guj is reflexive.

\subsection{Translators} The second extra plant,
besides of an adjoint \guj, must be another creator
$\star\in(0\mapsto 1)$ labelled by a star to be distinguished from the
first creator labelled by bullet $\bullet\in(0\mapsto 1).$ We are
showing the fragment of a garden expanded by two extra `dictionary'
plants, together with suitable new fundamental `dictionary' footpaths,
on Figure 16.

\guR 

We need also to give an excuse to the reader for denoting `translators'
paths by $d^*$ \& $\delta^*.$ These arrows correspond to the case when
a model for a plant \guj is a contra-variant functor, however, such
choice was not yet made. Moreover, such notation suggests that the
translators paths must be related to the Cartan fundamental paths $d$
\& $\delta.$ At this moment we prefer to consider such notation to be
nothing more than a convenient one. Possible not trivial
interconnections might appear after we introduce a Cartan-Leibniz \&
Leibniz-Cartan dictionary.

\section{Dictionary: from Cartan 2-cells to Leibniz 2-cells} In
order to translate the Cartan (`bimodule valued') fundamental path $d$
(as well as $\delta$ with source in a cobimodule) into the Leibniz-like
derivation, we need the translator paths $\delta^*$ and $d^*$ from
Figure 16.

\guS   

Figure 17 contain a subtle point, not shown for simplicity. The
evaluation \& coevaluation paths can hold for a unital binary plant
\guy\, only. Therefore a creator $\gcr\in(0\mapsto 1)$ on Figure 17
must be unit for evaluation \& coevaluation paths (as on Figures 14 \&
15). However \gcr\, needs not to be a unit for other paths on Figure 17.
Therefore in fact, what is missing on the top of Figure 17 is a path
`ch' between two different binary plants \guy\,, which must be labelled
(alternatively this can be traced to the change of the categories),
$$(\hbox{\guy for which \gcr\, is {\it not} unit})\quad
\stackrel{\hbox{ch}}{\lra}\quad(\hbox{unital \guy with unit \gcr\,}).$$

The Figure 17 is known in the literature as `the Cartan formula'.
For example, we can recover a left action $l$ and a left coaction $L$
from the Cartan garden, formally as follows
$$l\equiv\ev_l\circ\hbox{ch}\circ(\delta^*\times d),\qquad
L\equiv(d^*\times\delta)\circ\hbox{ch}^{-1}\circ\coev_r.$$
However, it is not clear here what could mean ch$^{-1}$?

Figure 17 is a categorification of a naive definition of a partial
`derivation' $\partial_\mu$ from the Cartan `derivation' (we do not yet
use any Leibniz-like axiom), viz.,
$$df\longmapsto\partial_\mu f=(df)\partial_\mu\quad(=\ev_l[(\delta^*
\partial_\mu)\otimes(df)]).$$
In the classical differential geometry a translator $\delta^*$ on Figure
16 is denoted by $i\equiv\delta^*,$ and moreover the following strict
3-cells-axioms holds (\ie $d^*$ \& $\delta^*$ are pseudo relations),
\begin{equation}d^*\circ\delta^*\approx\id\qquad\hbox{and}\qquad
\delta^*\circ d^*\approx\id.\label{cdg}\end{equation}
An unsolved problem in [Borowiec \& V\'azquez Couti\~no] can be
reformulated equivalently as the word-like problem of a `compatibility'
of the strict 3-cells (\ref{cdg}) with a Leibniz 3-cell `der'.

\section{Dictionary: from Leibniz 2-cells to Cartan bimodule}
In order to translate the Leibniz-like (right or left) (co)derivations
$\{r,l,R,L\}$ into the Cartan fundamental paths $d$ (target in a
bimodule $m_r,m_l$) \& $\delta$ (source in a cobimodule
$\triangle_r,\triangle_l$), we need the translator paths $d^*$ \&
$\delta^*$ from Figure 16.

\guT\smallskip  

In Figure 18, for abbreviation, we let a `unit' $u$ stand for the
composition of a left genuine unit$^{-1}$ for a bifunctor \guy\,
with a virtual unit $v$ for an algebra $m$ (virtual because we are not
yet sure that we need a compatibility axiom with a path $m$:
invertibility of the composed 2-cells $m\circ(v\times\id)$ \&
$m\circ(\id\times v)$), \ie we have the following abbreviations:
\begin{center}$u=\left(\right.$\bt virtual unit for\\an algebra $m$\et
$\left.\right)\circ\left(\right.$\bt a left unit$^{-1}$\\for a plant
\guy\et$\left.\right),$\\$\varepsilon=\left(\right.$\bt left unit for\\
a plant \guy\et$\left.\right)\circ\left(\right.$\bt a virtual counit\\
for a cogebra $\triangle$\et$\left.\right).$\end{center}
Figure 18 is the categorification of a well known naive `definition' of
the Cartan derivation $d\in\der(A,M)$ in terms of the given partial
derivations, viz.,
$$\partial_\mu\longmapsto df\equiv(dx^\mu)(\partial_\mu f).$$
The above `definition' is like a perpetum mobile because the
differential $df$ is `defined' again in terms of the differentials
$dx^\mu$ and a cicle is closed.

\subsection{Partial braid} On the way, from the given
2-cells in the garden for the partial derivation 3-cell, to the Cartan
club, there is the need of the direct constructions of an extended gebra
bi(co-)module 2-cells-paths $\{m_r,m_l,\triangle_r,\triangle_l\}$ from
the actions $\{r,l,R,L\}.$ The crucial step for this is the construction
of the derived 2-cells-paths of `partial braiding' of \gcr\, with $|$,
where a creator \gcr\, is {\it not} a unit for \guy, as shown here

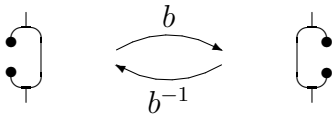
\guX  

One can show that $b$ can be given by the right action $r,$ and
$b^{-1}$ analogously by the right coaction $R.$ However, limited space
does not allow to present this construction here in detail. Also,
limited space does not allow to present the deeper consequences of our
assumptions.

\section{Axioms, word problem, coherence} Axioms generate binary
relations in the set of all expanded footpaths. The determination of
the category-congruence generated by axioms relating footpaths is almost
the Thue classical word problem [or the Birkhoff problem in universal
algebra: to determine all other `equations, or identities' as the
consequence of the given axioms], except that footpaths (our letters)
are not always composable. These two word problems, for plants and for
footpaths, are called jointly the coherence problem for a free club
[Kelly \& Laplaza 1980, \S 4 p. 198 and \S 10, p. 211].

Let $\approx$ be a binary relation among derived plants. As is usual
for the word problem, we use `relation' also for an element of an actual
binary equi\-valence relation, \ie also for an element $\in\;\approx,$
like the examples of generating relations on Figure 2. An equivalence
binary relation $\approx$ on derived plants is said to be a congruence
if $\approx$ is compatible with grafting, \ie if a grafting of a related
pair from $\approx$ gives again the related pair from $\approx.$
A category-congruence is generated from the given relations by means of
the Birkhoff's rules of derivations, see e.g. in [Graczy\'nska 1998, p.
13]. The determination of this congruence is the classical word problem of
Thue and this is a first part of the coherence problem for a club.
A free club factored by a congruence $\approx$ gives the quotient
category.

Axioms among plants (axioms generating the binary equivalence relations
$\approx$ on plants, axiom $\in\;\approx$), by definition, are allowed
within the same category only
$$\approx\quad\subset\quad(n\mapsto m)\times(n\mapsto m).$$
No relations are allowed between plants of the different arities. Axioms
among plants, as the three examples on Figure 2 within a category
$1\mapsto 1$, do express the properties of the involved plants in terms
of the derived plants.

An axiom is a particular case of the more general concept of a
modification, introduced by B\'enabou in 1967. According to this notion,
an axiom is the same as a strict modification. A weaker version is a
quasi-axiom = a quasi-modification, which is an invertible modification,
and still weaker is the so called `lax' version (most general) which is
said just to be a modification.

\subsection{Example of a Leibniz 3-cell} The Leibniz
relation is a 3-cell in an abelian operad generated by an alphabet of
footpaths. Because of very limited space, we will give an illustrative
example only. Let on Figure 12, $D$ be a composition of a fixed element
of $L$ (a creator in a model) with a left action $l.$ Then, a classical
example of a strict 3-cell for the set of 2-cells $\{D,m,\psi\}$ is
illustrated on Figure 19 in a model,
$m\equiv m_A,D\equiv D_A,\psi\equiv\psi_A,$

\guY\smallskip   

The first two terms (composed 2-cells) on Figure 19 are precisely the
Leibniz operations from Figure 12. However, the last term needs a
braiding of objects $L$ with $A,$ not included into our simplest
alphabet of footpaths in the present paper.

\begin{problem}[Problem] One of the open `higher order' problems in a
categorical \& algebraic analysis is the systematic investigations of a
weaker (lax) version of the Leibniz-like axioms (lax axioms for 2-cells,
either for partial derivatives as well as for the Cartan bimodule-valued
derivations), and the determination of the new strict axioms for the
Leibniz modifications, \ie the strict axioms for the Leibniz's 3-cells.
This would be the alternative, even the best, way to understand the
non-Leibniz component in [Przeworska-Rolewicz 1995].\end{problem}

\section{Functorial realization: functorial models} Eilenberg and
Mac Lane in 1945 introduced the notions of category, functor and natural
transformation of functors. In a functor category
$({\cat}_B)^{{\cat}_A},$ the $\set$-valued hombifunctor is denoted by
$\nat.$ Let $f,g$ be two functors from ${\cat}_A$ to ${\cat}_B,$
$f,g\in\obj\left\{({\cat}_B)^{{\cat}_A}\right\}.$
A natural transformation $t\in\nat(f,g)$ is an application
$t:\obj\cat_A\ra\cat_B(f\cdot,g\cdot),$ \ie a family of morphisms in
$\cat_B,$ such that every morphism $\phi\in\cat_A(\cdot,\cdot)$ gives
rise to a commutative diagram [Mac Lane 1963, 1965, 1971; Eilenberg and
Kelly 1966; Par\'e and Rom\'an 1998].

\begin{rem2} We wish that natural transformations
must be always closed under composition in
$$\nat\left[\nat(f,g)\times\nat(g,h),\nat(f,h)\right].$$
This implies that in the case of bivariant (\ie of mixed variances)
(multi)functors the definition of a natural transformation, as given by
Eilenberg \& Kelly in 1966, and by Kelly in [1972, \S 4 pp.93-94], must
be modified. A `generalized' natural transformation is known under the
names: dinatural transformation, or Barr dinatural transformation
[Par\'e \& Rom\'an 1998].\end{rem2}

A bicategory [B\'enabou 1967] is the same as a not strict 2-category
(also called pseudo or more general lax). A typical example of a
2-category is a 1-category, or, in the case of many sorts, the
collection of 1-categories, {\it with a structure,} according to the
following definition.

\begin{def4}[Kelly 1972, \S 3, p. 116]. A {\it structure}\/
on a category consists of a set of functors, and of various natural
transformations, subject to equational axioms (like, for example, the
Leibniz-like axiom).\end{def4}

Lawvere in his Thesis in 1963 gave a program of a categorification: a
program of replacing standard `theories', \ie a language with rules of
deduction and axioms, by n-categories with structures, `models' by
n-functors and morphisms by natural transformations.
A theory, a mathematical theory, such as the algebraic theory [Lawvere
1963] based on Birkhoff's equational class, or in particular an
algebraic analysis (based on the Leibniz-like axiom), consists of
an n-category (or n-categories = sorts in computer science) with
various amounts of structure imposed on them. For the most elementary
introduction, see [Lawvere \& Schanuel 1997]. For a most recent account
we can refer to [Baez \& Dolan 1998, Batanin 1998].

An n-functor from an algebraic n-category {\bf N} with grafting, into
an n-category {\bf Cat} of small functor categories (with grafting), is
said to be the functorial model of {\bf N} (a realization or an
{\bf N}-algebra). Such n-functor is a model of plants as 1-functors, and
a model of footpaths as natural transformations of 1-functors.

\section{Conclusions and some directions for further studies}
The authors see an algebraic analysis as just another instance of an
equational theory (an instance of the Birkhoff universal algebra) with
a Leibniz-like axiom as a strict or a weak 3-cell. In this paper we do
propose to apply Lawvere's categorification program to analysis. Very
limited space does not allow to discuss here neither the specific
Leibniz 3-categories with weak 3-cells nor many other important issues.
Therefore our purpose was limited to an elementary introduction into
the basic notions of the categorical analysis in terms of the general
theory of n-categories. There are many important directions which
deserve further detailed studies, besides of these mentioned in the
main body of this paper. A logical order of these directions, from the
most fundamental to the less relevant, is contrary to the historical
developments.

\subsection{Extension of gebra} We believe that the
most fundamental for any analysis (algebraic, categorical, braided) is
an extension of gebra \ie an extension of 2-cells in terms of 3-cells
denoted below by `$\longmapsto$', [Eilenberg 1948, Gugenheim 1962,
Kelly 1972, p. 94, Cuntz \& Quillen 1995], viz.,
$$m\longmapsto m_l\,\&\,m_r\qquad\hbox{or/and}\qquad \triangle
\longmapsto\triangle_l\,\&\,\triangle_r.$$ The {\it extension} 3-cells
are said to be axioms if they are strict, and have nothing to do with a
grafting dependent {\it expansion} of 2-cells or with a B\'enabou's
pasting of 2-cells. The strict Eilenberg-like extension 3-cells are
known in the non-commutative differential geometry under the name of
gebra bimodule, bicomodule for cogebra, double dimodule [Pareigis 1996],
quadruple comodule [Borowiec \& V\'azquez Couti\~no], etc. A gebra
extension must precede algebraic analysis understood as a Leibniz-like
3-cell. We see that the logical order of differential calculus is
as follows: firstly the 3-cells of an Eilenberg-like gebra extension,
and after a Leibniz-like 3-cell.

The extension 3-cells need an abelian category, \ie a pair of binary
plants (bifunctors) and moreover a distributive law between them [Beck
1969; Kelly 1972, p. 94].

\subsection{A Leibniz-like axiom as a 3-cell}
Classical differential geometry deals with the Cartan derivationm from
an algebra $m$ to an $m$-bimodule $m_l$ \& $m_r$ of the differential
one-forms, $d\in\der(m,m_l\,\&\,m_r).$ Therefore $\der$ is interpreted
usually as a bifunctor
$$\hbox{(algebra,bimodule)}\quad\longmapsto\quad\set.$$
We believe that it would be desirable to identify der also as a strict
3-cell in an appropriate Leibniz 3-category and to investigate a weak
version in terms of the B\'enabou modification.

\subsection{Braided analysis} The classical work on a
braided analysis was done by Woronowicz [1989], who located partial
braided derivations $\der(m,\sigma,\ldots)$ within a braided Lie
algebra. Majid in [1993, 1995 \S 10.4] introduced another braided
version of the Leibniz rule, that is, another braided derivation. From
the point of view of the categorical analysis, which we do propose here,
Majid's braided analysis seems to be not so much an important example.
A slightly different approach was proposed in [Oziewicz, Paal \&
R\'o\.za\'nski 1995]. The main unsolved problem here is to locate
partial derivations (braided or not braided), as a 3-cell
$\der(m)\equiv\der(m,m\,\&\,m),$ within an extra structure of a braided
Lie algebra which we do propose to identify as some strict 4-cell, \ie
as an axiom on 3-cells. In particular, it would be very interesting to
identify possible strict 4-cells for the Borowiec partial
(co)derivations [Borowiec 1996, 1997; Borowiec \& V\'azquez Couti\~no].

\end{document}